\documentclass[11pt]{article}

\usepackage{amsfonts,amssymb,amsmath,amsthm,epsfig,euscript}

\newtheorem{theorem}{Theorem}
\newtheorem{lemma}[theorem]{Lemma}

\long\def\symbolfootnote[#1]#2{\begingroup
\def\thefootnote{\fnsymbol{footnote}}\footnote[#1]{#2}\endgroup}

\newcommand{\la}{\lambda}
\newcommand{\La}{\Lambda}

\newcommand{\des}[1]{\mathrm{des}(#1)}
\newcommand{\cdes}[1]{\mathrm{cdes}(#1)}
\newcommand{\cyc}[1]{\mathrm{cyc}(#1)}
\newcommand{\LRmin}[1]{\mathrm{LRmin}(#1)}

\newcommand{\nth}[1][n]{{#1}^{\mathrm{th}}}

\newcommand{\sg}{\sigma}

\newcommand{\cref}[1]{Corollary \ref{corollary:#1}}

\newcommand{\red}[1]{\mathrm{red}(#1)}

\newcommand{\tmch}[1]{\text{$\tau$-$\mathrm{mch}$}(#1)}
\newcommand{\ctmch}[1]{\text{$c$-$\tau$-$\mathrm{mch}$}(#1)}

\title{A reciprocity method for computing generating functions over the set of permutations with no consecutive occurrence of $\tau$}

\author{
Miles Eli Jones \\
\small Department of Mathematics\\[-0.8ex]
\small University of California, San Diego\\[-0.8ex]
\small La Jolla, CA 92093-0112. USA\\[-0.8ex]
\small \texttt{mej005@math.ucsd.edu}
\and
\and
Jeffrey B. Remmel \\
\small Department of Mathematics\\[-0.8ex]
\small University of California, San Diego\\[-0.8ex]
\small La Jolla, CA 92093-0112. USA\\[-0.8ex]
\small \texttt{remmel@math.ucsd.edu}
\and
}

\date{\small Submitted: Date 1;  Accepted: Date 2;
 Published: Date 3.\\
\small MR Subject Classifications: 05A15, 05E05 \\
keywords: permutation, pattern match, descent, left to right minimum, symmetric polynomial, exponential generating function}

\begin{document}

\maketitle
 
\begin{abstract} In this paper, we introduce a new method for 
computing generating functions with respect to the number of descents 
and left-to-right minima over the set of permutations which have 
no consecutive occurrence of $\tau$ where $\tau$ starts with 1. 
In particular, we study the generating function 
$\sum_{n \geq 0} \frac{t^n}{n!} 
\sum_{\sg \in \mathcal{NM}_n(1324 \ldots p)}x^{\LRmin{\sg}}y^{1+\des{\sg}}$ where $p \geq 4$, 
$\mathcal{NM}_n(1324 \ldots p)$  
is the set of permutations $\sg$ in the symmetric group $S_n$ which has 
no consecutive occurrences of $1324 \ldots p$, $\des{\sg}$ is the number 
of descents of $\sg$ and $\LRmin{\sg}$ is the number of left-to-right 
minima of $\sg$. We show that for any $p \geq 4$, this generating function  
is of the form $\left( \frac{1}{U(t,y)}\right)^x$ where 
$U(t,y) = \sum_{n\geq 0}U_n(y) \frac{t^n}{n!}$ and the coefficients 
$U_n(y)$ satisfy some simple recursions depending on $p$.  
As an application of our results, we compute explicit generating 
functions for the number of permutations of $S_n$ that have 
no consecutive occurrences of the pattern $1324 \ldots p$ and 
have exactly $k$ descents for $k=1,2$. 
\end{abstract}

\section{Introduction}

Given a sequence \begin{math}\sg = \sg_1 \ldots \sg_n\end{math} of distinct integers,
let \begin{math}\red{\sg}\end{math} be the permutation found by replacing the
\begin{math}i^{\textrm{th}}\end{math} largest integer that appears in \begin{math}\sg\end{math} by \begin{math}i\end{math}.  For
example, if \begin{math}\sg = 2~7~5~4\end{math}, then \begin{math}\red{\sg} = 1~4~3~2\end{math}.  Given a
permutation \begin{math}\tau=\tau_1 \ldots \tau_j\end{math} in the symmetric group \begin{math}S_j\end{math}, we say a
permutation \begin{math}\sg = \sg_1 \ldots \sg_n \in S_n\end{math} has  a \emph{
\begin{math}\tau\end{math}-match starting at position \begin{math}i\end{math}} provided \begin{math}\red{\sg_i \ldots \sg_{i+j-1}}
= \tau\end{math}.  Let \begin{math}\tmch{\sg}\end{math} be the number of \begin{math}\tau\end{math}-matches in the
permutation \begin{math}\sg\end{math}.  Given a permutation \begin{math}\sg = \sg_1 \ldots \sg_n \in S_n\end{math}, we let 
\begin{math}\des{\sg} = |\{i: \sg_i > \sg_{i+1}\}|\end{math}.  We say that 
\begin{math}\sg_j\end{math} is a \emph{left-to-right minimum} of \begin{math}\sg\end{math} if \begin{math}\sg_j < \sg_i\end{math} for 
all \begin{math}i<j\end{math}. We let \begin{math}\LRmin{\sg}\end{math} denote the number of 
left-to-right minima of \begin{math}\sg\end{math}. 
Let \begin{math}\mathcal{NM}_n(\tau)\end{math} denote the set of permutations in \begin{math}S_n\end{math} with no \begin{math}\tau\end{math}-matches and 
let 
\begin{equation}
NM_{\tau,n}(x,y) = 
\sum_{\sg \in \mathcal{NM}_n(\tau)} x^{\LRmin{\sg}} y^{1+\des{\sg}}.
\end{equation}

The main goal of this paper is to give a new method  to 
compute generating 
functions of the form 
\begin{equation}
NM_{\tau}(t,x,y) = \sum_{n \geq 0} \frac{t^n}{n!} NM_{\tau,n}(x,y) 
\end{equation}
where where $\tau$ is a permutation in $S_j$ which starts with 1. 
This method was first outlined in the \cite{JR11} where 
the authors considered the generating function $NM_{1423}(t,x,y)$. 
In this paper, we shall develop a more general approach that applies 
to any $\tau$ that starts with 1 and has one descent. 
Our method does not compute \begin{math}NM_{\tau}(t,x,y)\end{math} directly.  
Instead, we assume that     
\begin{equation}\label{eq:I1}
NM_{\tau}(t,x,y) = \left(\frac{1}{U_\tau(t,y)}\right)^x
\ \mbox{where} \ 
U_\tau(t,y) = 1 + \sum_{n\geq 1} U_{\tau,n}(y) \frac{t^n}{n!}.
\end{equation}
Thus 
\begin{equation}\label{eq:I3}
U_{\tau}(t,y) = \frac{1}{1+\sum_{n \geq 1} NM_{\tau,n}(1,y) \frac{t^n}{n!}}
\end{equation}
Our method gives a combinatorial 
interpretation the right-hand side of (\ref{eq:I3}) and then uses  
that combinatorial interpretation to develop simple recursions on 
the coefficients $U_{\tau,n}(y)$.

The first question to ask is why can we assume that $NM_{\tau}(t,x,y)$ 
can be expressed in the form 
$$NM_{\tau}(t,x,y) = \left(\frac{1}{U_\tau(t,y)}\right)^x$$
when $\tau$ starts with 1.  The fact that this is assumption is 
justified 
follows from previous work of the authors \cite{JR} where they 
introduced the general study of patterns in 
the cycle structure of permutations. 
That is, suppose that 
\begin{math}\tau=\tau_1 \ldots \tau_j\end{math} is a permutation in \begin{math}S_j\end{math} and 
\begin{math}\sg\end{math} is a permutation in \begin{math}S_n\end{math} with \begin{math}k\end{math} cycles \begin{math}C_1 \ldots C_k\end{math}. 
We shall always write cycles in the form 
\begin{math}C_i =(c_{0,i}, \ldots, c_{p_i-1,i})\end{math}  where \begin{math}c_{0,i}\end{math} is the smallest 
element in \begin{math}C_i\end{math} and \begin{math}p_i\end{math} is the length of \begin{math}C_i\end{math}. Then we arrange 
the cycles by decreasing smallest elements. That is, we arrange 
the  cycles of \begin{math}\sg\end{math} so that \begin{math}c_{0,1} > \cdots > c_{0,k}\end{math}. Then 
we say that \begin{math}\sg\end{math} has a {\em cycle-\begin{math}\tau\end{math}-match} (\begin{math}c\end{math}-\begin{math}\tau\end{math}-mch)
if there is an \begin{math}i\end{math} such that \begin{math}C_i =(c_{0,i}, \ldots, c_{p_i-1,i})\end{math} where 
\begin{math}p_i \geq j\end{math} and an \begin{math}r\end{math} such that 
\begin{math}\red{c_{r,i} c_{r+1,i} \ldots c_{r+j-1,i}} = 
\tau\end{math} where we take indices of the form \begin{math}r+s\end{math} modulo \begin{math}p_i\end{math}.   
Let \begin{math}\ctmch{\sg}\end{math} be the number of cycle-\begin{math}\tau\end{math}-matches in the
permutation \begin{math}\sg\end{math}.  For example, 
if \begin{math}\tau =2~1~3\end{math} and \begin{math}\sg = (4,7,5,8,6)(2,3)(1,10,9)\end{math}, then 
\begin{math}9~1~10\end{math} is a cycle-\begin{math}\tau\end{math}-match in the third cycle and 
\begin{math}7~5~8\end{math} and \begin{math}6~4~7\end{math} are cycle-\begin{math}\tau\end{math}-matches in the first cycle so 
that \begin{math}\ctmch{\sg} =3\end{math}.  

Given a cycle 
\begin{math}C=(c_0, \ldots, c_{p-1})\end{math} where \begin{math}c_0\end{math} is the smallest element in 
the cycle, we let \begin{math}\cdes{C} = 1+ \des{c_0 \ldots c_{p-1}}\end{math}. Thus 
\begin{math}\cdes{C}\end{math} counts the number of descent pairs as 
we traverse 
once around the cycle because the extra  \begin{math}1\end{math} counts 
the descent pair \begin{math}c_{p-1}>c_0\end{math}. For example if \begin{math}C = (1,5,3,7,2)\end{math}, then 
\begin{math}\cdes{C} = 3\end{math} which counts the descent pairs \begin{math}53\end{math}, \begin{math}72\end{math}, and \begin{math}21\end{math} as 
we traverse once around \begin{math}C\end{math}.  By convention, if 
\begin{math}C =(c_0)\end{math} is a one-cycle, we let \begin{math}\cdes{C} = 1\end{math}. If \begin{math}\sg\end{math} is a permutation in \begin{math}S_n\end{math} with \begin{math}k\end{math} cycles \begin{math}C_1 \ldots C_k\end{math}, then 
we define \begin{math}\cdes{\sg} = \sum_{i=1}^k \cdes{C_i}\end{math}.  We let 
\begin{math}\cyc{\sg}\end{math} denote the number of cycles of \begin{math}\sg\end{math}.

In \cite{JR}, 
Jones and Remmel studied generating functions of the form 
\begin{equation*}
NCM_{\tau}(t,x,y) =  1 + \sum_{n \geq 1}  \frac{t^n}{n!} 
\sum_{\sg \in \mathcal{NCM}_n(\tau)} x^{\cyc{\sg}} y^{\cdes{\sg}} 
\end{equation*}
where  
\begin{math}\mathcal{NCM}_n(\tau)\end{math} is 
the set of permutations \begin{math}\sg \in S_n\end{math} which have 
no cycle-\begin{math}\tau\end{math}-matches. 
The basic approach used in that paper was to use 
theory of exponential structures to reduce the problem 
of computing \begin{math}NCM_{\tau}(t,x,y)\end{math} 
to the problem of computing similar generating functions for \begin{math}n\end{math}-cycles. 
That is, let \begin{math}\mathcal{NCM}_{n,k}(\tau)\end{math} be the set  of permutations \begin{math}\sg\end{math} 
of \begin{math}S_n\end{math} with \begin{math}k\end{math} cycles such that \begin{math}\sg\end{math} has no cycle-\begin{math}\tau\end{math}-matches 
and \begin{math}\mathcal{L}_m^{ncm}(\tau)\end{math} denote the set  of \begin{math}m\end{math}-cycles \begin{math}\gamma\end{math} 
in \begin{math}S_m\end{math} such \begin{math}\gamma\end{math} has no cycle-\begin{math}\tau\end{math}-matches. 
The following theorem follows easily from the theory of exponential 
structures as is described in \cite{Stan}, for example.                                    

\begin{theorem}\label{Fundrefined} 
\begin{equation}\label{expncmUtxy}
NCM_{\tau} (t,x,y) = 1 + \sum_{n \geq 1} \frac{t^n}{n!} \sum_{k=1}^n 
x^k \sum_{\sg \in \mathcal{NCM}_{n,k}(\tau)} y^{\cdes{\sg}}= 
e^{x \sum_{m \geq 1}  \frac{t^m}{m!} 
\sum_{C \in \mathcal{L}^{ncm}_m(\tau)} y^{\cdes{C}}}.
\end{equation}
\end{theorem}

It turns out that if \begin{math}\tau \in S_j\end{math} is a permutation that 
starts with 1, then we can reduce the problem of finding 
\begin{math}NCM_{\tau}(t,x,y)\end{math} to the usual problem of finding 
the generating function of permutations that have no \begin{math}\tau\end{math}-matches. 
Let \begin{math}\bar{\sg}\end{math} be the permutation that arises from 
\begin{math}C_1 \cdots C_k\end{math} by erasing all the parentheses and commas. 
For example, if 
\begin{math}\sg = (7,10,9,11)\ (4,8,6) \ (1,5,3,2)\end{math}, then 
\begin{math}\bar{\sg} = 7~10~9~11~4~8~6~1~5~3~2\end{math}. It is easy to see that 
the minimal elements of the cycles correspond to left-to-right minima in 
\begin{math}\bar{\sg}\end{math}.  It is also easy to see 
that under our bijection \begin{math}\sg \rightarrow \bar{\sg}\end{math}, 
\begin{math}\cdes{\sg} = \des{\bar{\sg}}+1\end{math} since every left-to-right minima other than the first element of $\bar{\sg}$ 
is part of a descent pair in \begin{math}\bar{\sg}\end{math}. For example, if 
\begin{math}\sg = (7,10,9,11)\ (4,8,6) \ (1,5,3,2)\end{math} so that 
\begin{math}\bar{\sg} = 7~10~9~11~4~8~6~1~5~3~2\end{math}, 
\begin{math}\cdes{(7,10,9,11)} = 2\end{math}, \begin{math}\cdes{(4,8,6)} =2\end{math}, and 
\begin{math}\cdes{(1,5,3,2)} =3\end{math} so that \begin{math}\cdes{\sg} = 2+2+3 =7\end{math} while 
\begin{math}\des{\bar{\sg}} =  6\end{math}.  This given, Jones and Remmel 
\cite{JR} proved that if \begin{math}\tau \in S_j\end{math} and \begin{math}\tau\end{math} starts with 1, then for any \begin{math}\sg \in S_n\end{math},  
(1) \begin{math}\sg\end{math} has \begin{math}k\end{math} cycles if and only if \begin{math}\bar{\sg}\end{math} has \begin{math}k\end{math} left-to-right
minima, (2) \begin{math}\cdes{\sg} = 1+\des{\bar{\sg}}\end{math}, and 
(3) \begin{math}\sg\end{math} has no cycle-\begin{math}\tau\end{math}-matches if and 
only if \begin{math}\bar{\sg}\end{math} has no \begin{math}\tau\end{math}-matches.  Thus they proved the following theorem. 
\begin{theorem} \label{thm:tau1}
Suppose that $\tau = \tau_1 \ldots \tau_j \in S_j$ and 
$\tau_1 =1$. Then 
\begin{equation}
NCM_{\tau}(t,x,y) = NM_{\tau}(t,x,y). 
\end{equation}
\end{theorem}

It follows from Theorems \ref{Fundrefined} and \ref{thm:tau1} 
that if $\tau \in S_j$ and $\tau$ starts with 1, then 
\begin{equation}
NM_{\tau}(t,x,y) =F(t,y)^x
\end{equation}
 for some function $F(t,y)$. Thus our assumption that 
\begin{equation}
NM_{\tau}(t,x,y) =\left(\frac{1}{U_\tau(t,y)}\right)^x
\end{equation}
is fully justified in the case when $\tau$ starts with 1. 
We should note that if a permutation \begin{math}\tau\end{math} does not start 
with 1, then it may be that case that \begin{math}|\mathcal{NM}_n(\tau)| \neq 
|\mathcal{NCM}_n(\tau)|\end{math}. For example, Jones and Remmel \cite{JR} computed that 
\begin{math}|\mathcal{NCM}_7(3142)| = 4236\end{math} and \begin{math}
|\mathcal{NM}_7(3142)| = 4237\end{math}.

Jones and Remmel \cite{JR} were able 
compute functions of the form \begin{math}NCM_{\tau}(t,x,y)\end{math} when \begin{math}\tau\end{math} starts with 
1 by combinatorially  
proving certain recursions for $\sum_{C \in \mathcal{L}^{ncm}_m(\tau)} y^{\cdes{C}}$ which led 
to certain sets of differential equations satisfied 
by  \begin{math}NCM_{\tau}(t,x,y)\end{math}. 
For example, using such methods, 
they were able to prove the following two theorems.

\begin{theorem}\label{thm:1-2} 
Let \begin{math}\tau = \tau_1 \ldots \tau_j \in S_j\end{math} where 
\begin{math}j \geq 3\end{math} and \begin{math}\tau_1 =1\end{math} and \begin{math}\tau_j =2\end{math}.
Then 
\begin{equation}\label{1alpha2}
NCM_{\tau}(t,x,y) = NM_{\tau}(t,x,y) =\frac{1}{\left(1 - \int_0^t e^{(y-1)s- \frac{y^{\des{\tau}}s^{j-1}}{(j-1)!}}ds\right)^x}
\end{equation}
\end{theorem}

\begin{theorem}\label{thm:1-j-1-sg-j} Suppose that \begin{math}\tau = 12\ldots (j-1)(\gamma) j\end{math} where 
\begin{math}\gamma\end{math} is a permutation of \begin{math}j+1, \ldots , j+p\end{math} and \begin{math}j \geq 3\end{math}. 
Then 
\begin{displaymath}NCM_\tau(t,x,y) =  NM_\tau(t,x,y) = \frac{1}{(U_\tau(t,y))^x}\end{displaymath} 
where 
\begin{math}\displaystyle U_\tau(t,y) = 1+\sum_{n \geq 1} U_{\tau,n}(y) \frac{t^n}{n!}\end{math}, $U_1(y) =-y$, 
and for all \begin{math}n \geq 2\end{math}, 
\begin{equation}\label{Urec}
U_{\tau,n}(y) = (1-y)U_{\tau,n-1}(y) - y^{\des{\tau}}\binom{n-j}{p} 
U_{\tau,n-p-j+1}(y).
\end{equation}
\end{theorem}

The next step in our approach is to use the homomorphism method 
to give a combinatorial interpretation to 
\begin{equation}
U_\tau(t,y) = \frac{1}{NM_\tau(t,1,y)}.
\end{equation}
That is, Remmel and various  coauthors \cite{Bec4,Lan2,MenRem1,MenRem2,MenRem3,MRR,Remriehl} developed a method called 
the homomorphism method to show that many generating 
functions involving permutation statistics can be derived by 
applying a homomorphism defined on the 
ring of symmetric functions \begin{math}\Lambda\end{math}  
in infinitely many variables \begin{math}x_1,x_2, \ldots \end{math} 
to simple symmetric function identities such as 
\begin{equation}\label{conclusion2}
H(t) = 1/E(-t)
\end{equation}
where 
\begin{equation}\label{genfns}
H(t) = \sum_{n\geq 0} h_n t^n = \prod_{i\geq 1} \frac{1}{1-x_it} \ \mbox{and} \ E(t) = \sum_{n\geq 0} e_n t^n = \prod_{i\geq 1} 1+x_it
\end{equation}
are the generating functions of the homogeneous symmetric functions  
\begin{math}h_n\end{math} and the elementary symmetric functions  
\begin{math}e_n\end{math} in infinitely many variables \begin{math}x_1,x_2, \ldots \end{math}. 
Now suppose that we define a homomorphism \begin{math}\theta\end{math} on 
\begin{math}\Lambda\end{math} by setting 
\begin{displaymath}\theta(e_n) = \frac{(-1)^n}{n!} NM_{\tau,n}(1,y).\end{displaymath}
Then 
\begin{displaymath}\theta(E(-t)) = {\sum_{n\geq 0} NM_{\tau,n}(1,y) \frac{t^n}{n!}} = \frac{1}{U_\tau(t,y)}.\end{displaymath}
Thus \begin{math}\theta(H(t))\end{math} should equal \begin{math}U_\tau(t,y)\end{math}.  We shall then show 
how to use the combinatorial methods associated with 
the homomorphism method to develop recursions for 
the coefficient of \begin{math}U_\tau(t,y)\end{math} similar to those  
in Theorem \ref{thm:1-j-1-sg-j}.

For example, in this paper, we shall study the generating 
functions for permutations $\tau$ of the form 
$\tau = 1324\ldots p$ where $p \geq 4$. That is, $\tau$ arises 
from the identity 
by interchanging 2 and 3. We shall show that $U_{1324,1}(y)=-y$ and 
for $n \geq 2$, 
\begin{equation}\label{eq:1324}
U_{1324,n}(y) = (1-y)U_{1324,n-1}(y)+ \sum_{k=2}^{\lfloor n/2 \rfloor} 
(-y)^{k-1} C_{k-1} U_{1324,n-2k+1}(y)
\end{equation}
where \begin{math}C_k\end{math} is the \begin{math}k\end{math}-th 
Catalan number.  For any $p \geq 5$, we shall show 
that $U_{1324 \ldots p,n}(y) =-y$ and for $n \geq 2$, 
\begin{equation}\label{eq:1324p}
U_{1324\ldots p,n}(y)=(1-y)U_{1324\ldots p,n-1}(y)+\sum_{k=1}^{\lfloor\frac{n-2}{p-2}\rfloor}(-y)^{k}U_{1324\ldots p,n-(k(p-2)+1)}(y).
\end{equation}

The outline of this paper is as follows. 
In Section 2, we shall briefly recall the background in the 
theory of symmetric functions that we will need for our proofs. 
Then in Section 3, we shall illustrate our method by 
proving (\ref{eq:1324}) and (\ref{eq:1324p}). Then in Section 4, 
we shall show how the results of Section 3 allow us to 
compute explicit generating 
functions for the number of permutations of $S_n$ that have 
no consecutive occurrences of the pattern $1324 \ldots p$ and 
exactly $k$ descents for $p \geq 4$ and $k=1,2$. Finally, in Section 5, 
we shall discuss our conclusions as well as some other results that follow by using the same method.

\section{Symmetric functions.}
In this section, we give the necessary background on 
symmetric functions that will be needed for our proofs.

Given a partition \begin{math}\la = (\la_1, \ldots \la_\ell)\end{math} where 
\begin{math}0 < \la_1 \leq \cdots \leq \la_\ell\end{math}, we let 
\begin{math}\ell(\lambda)\end{math} be the number of nonzero integers in \begin{math}\la\end{math}.  If the sum of
these integers is equal to \begin{math}n\end{math}, then we say \begin{math}\la\end{math} is a partition of \begin{math}n\end{math} and 
write \begin{math}\la \vdash n\end{math}.  

Let \begin{math}\Lambda\end{math} denote the ring of symmetric functions in infinitely  
many variables \begin{math}x_1,x_2, \ldots \end{math}. The \begin{math}\nth\end{math} elementary symmetric function \begin{math}e_n = e_n(x_1,x_2, \ldots )\end{math}  and \begin{math}\nth\end{math} homogeneous 
 symmetric function \begin{math}h_n = h_n(x_1,x_2, \ldots )\end{math} are defined by the generating functions given in (\ref{genfns}). 
For any partition \begin{math}\la = (\la_1,\dots,\la_\ell)\end{math}, let \begin{math}e_\la = e_{\la_1} 
\cdots e_{\la_\ell}\end{math} and \begin{math}h_\la = h_{\la_1} 
\cdots h_{\la_\ell}\end{math}.  It is well known that \begin{math}\{e_\la : \text{\begin{math}\la\end{math} is a 
partition}\}\end{math} is a basis for \begin{math}\La\end{math}.  In particular, \begin{math}e_0,e_1, \ldots \end{math} is 
an algebraically independent set of generators for \begin{math}\La\end{math} and, hence, 
a ring homomorphism \begin{math}\theta\end{math} on \begin{math}\Lambda\end{math} can be defined  
by simply specifying \begin{math}\theta(e_n)\end{math} for all \begin{math}n\end{math}.

\begin{figure}[htbp]
  \begin{center}
    \includegraphics[width=0.6\textwidth]{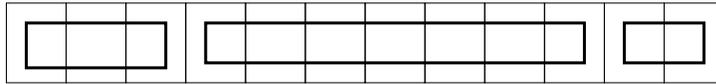}
    \caption{A $(2,3,7)$-brick tabloid of shape $(12)$.}
    \label{fig:DIMfig1}
  \end{center}
\end{figure}

A key element of our proofs is the combinatorial 
description of the coefficients of the expansion of 
\begin{math}h_n\end{math} in terms of the elementary symmetric functions 
\begin{math}e_\lambda\end{math} given by E\u{g}ecio\u{g}lu and Remmel in \cite{Eg1}.
They defined a \begin{math}\la\end{math}-brick tabloid of shape 
\begin{math}(n)\end{math} to be a rectangle of height \begin{math}1\end{math} and length \begin{math}n\end{math} which is covered by  ``bricks'' 
of  lengths found in the partition \begin{math}\la\end{math}
in such a way that no two bricks overlap. For example, Figure 
\ref{fig:DIMfig1} shows one brick \begin{math}(2,3,7)\end{math}-tabloid of shape \begin{math}(12)\end{math}.

Let \begin{math}\mathcal{B}_{\la,n}\end{math} denote the set of \begin{math}\la\end{math}-brick tabloids 
of shape \begin{math}(n)\end{math} and let \begin{math}B_{\la,n}\end{math} be the number of \begin{math}\la\end{math}-brick 
tabloids of shape \begin{math}(n)\end{math}.  If \begin{math}B \in \mathcal{B}_{\la,n}\end{math} we 
will write \begin{math}B =(b_1, \ldots, b_{\ell(\la)})\end{math} if the lengths of 
the bricks in \begin{math}B\end{math}, reading from left to right, are 
\begin{math}b_1, \ldots, b_{\ell(\la)}\end{math}. 
Through simple recursions, E\u{g}ecio\u{g}lu and Remmel  \cite{Eg1} proved 
that 
\begin{equation}\label{htoe}
h_n = \sum_{\la \vdash n} (-1)^{n - \ell(\la)} B_{\la,n} e_\la.
\end{equation}

\section{Computing $U_{1324 \ldots p,n}(y)$.}

The main goal of this section is to prove (\ref{eq:1324}) and 
(\ref{eq:1324p}). We shall start by proving (\ref{eq:1324p}).

Suppose that $\tau \in S_j$ which starts with 1 and $\des{\tau} =1$. 
Our first step is to give a combinatorial interpretation to 
\begin{equation}\label{eq:basic}
U_{\tau}(t,y)=\frac{1}{NM_{\tau}(t,1,y)}  = \frac{1}{1+ \sum_{n \geq 1} 
\frac{t^n}{n!} NM_{\tau,n}(1,y)}
\end{equation}
where $NM_{\tau,n}(1,y) = \sum_{\sg \in \mathcal{NM}_n(\tau)} y^{1+\des{\sg}}$.

To this end, we define a ring homomorphism $\theta_\tau$ on the ring of symmetric functions 
\begin{math}\Lambda\end{math} by setting 
\begin{math}\theta_{\tau}(e_0) =1\end{math} and  
\begin{equation}\label{theta}
\theta_{\tau}(e_n) = \frac{(-1)^n}{n!} NM_{\tau,n}(1,y) \ \mbox{for} \  n \geq 1.
\end{equation}
It follows that 
\begin{eqnarray*}
\theta_{\tau}(H(t)) &=& \sum_{n \geq 0} \theta_{\tau}(h_n)t^n  = \frac{1}{\theta_{\tau}(E(-t))} = \frac{1}{1 + \sum_{n \geq 1} (-t)^n \theta_{\tau}(e_n)} \\
&=& \frac{1}{1 + \sum_{n \geq 1} \frac{t^n}{n!} NM_{\tau,n}(1,y)}
\end{eqnarray*}
which is what we want to compute. 

\begin{figure}[htbp]
  \begin{center}
    \includegraphics[width=0.6\textwidth]{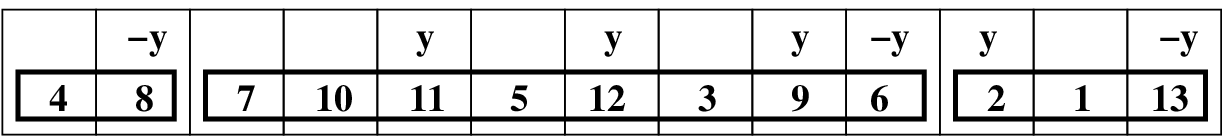}
    \caption{An element of $\mathcal{O}_{13245,13}$.}
    \label{fig:On}
  \end{center}
\end{figure}

By (\ref{htoe}), we have  that 
\begin{eqnarray}\label{eq:basic1}
n! \theta_{\tau}(h_n) &=& n! \sum_{\mu \vdash n} (-1)^{n-\ell(\mu)} 
B_{\mu,n} \theta_{\tau}(e_\mu) \nonumber \\
&=& n! \sum_{\mu \vdash n} (-1)^{n-\ell(\mu)} \sum_{(b_1, \ldots, 
b_{\ell(\mu)}) \in \mathcal{B}_{\mu,n}} \prod_{i=1}^{\ell(\mu)}  
\frac{(-1)^{b_i}}{b_i!} NM_{\tau,b_i}(1,y) \nonumber \\
&=& \sum_{\mu \vdash n} (-1)^{\ell(\mu)} \sum_{(b_1, \ldots, b_{\ell(\mu)}) \in \mathcal{B}_{\mu,n}} \binom{n}{b_1, \ldots, b_{\ell(\mu)}}
\prod_{i=1}^{\ell(\mu)}  NM_{\tau,b_i}(1,y).
\end{eqnarray}
Our next goal is to give a combinatorial interpretation to the 
right-hand side of (\ref{eq:basic1}).  If we are given 
a brick tabloid \begin{math}B=(b_1, \ldots, b_{\ell(\mu)})\end{math}, 
then we can interpret the multinomial coefficient \begin{math}\binom{n}{b_1, \ldots, b_\mu}\end{math} as all ways to assign sets \begin{math}S_1, \ldots , S_{\ell(\mu)}\end{math} to 
the bricks  of \begin{math}B\end{math} in such a way that \begin{math}|S_i| =b_i\end{math} for \begin{math}i=1, \ldots, \ell(\mu)\end{math} and the sets \begin{math}S_1, \ldots , S_{\ell(\mu)}\end{math} form a set partition 
of \begin{math}\{1, \ldots, n\}\end{math}.  Next for each brick \begin{math}b_i\end{math}, we use 
the factor \begin{displaymath}NM_{\tau,b_i}(1,y) = \sum_{\sg \in S_{b_i}, \tau\text{-mch}(\sg) =0} y^{\des{\sg}+1}\end{displaymath} to pick a rearrangement \begin{math}\sg^{(i)}\end{math} of \begin{math}S_i\end{math} which 
has no $\tau$-matches to put in cells of \begin{math}b_i\end{math} and then we place 
a label of \begin{math}y\end{math} on each cell that starts a descent in \begin{math}\sg^{(i)}\end{math} plus 
a label of \begin{math}y\end{math} on the last cell of \begin{math}b_i\end{math}. Finally, we use 
the term \begin{math}(-1)^{\ell(\mu)}\end{math} to turn each label \begin{math}y\end{math} at the end of 
brick to a \begin{math}-y\end{math}. We let \begin{math}\mathcal{O}_{\tau,n}\end{math} denote the set of 
all objects created in this way. For each element \begin{math}O \in \mathcal{O}_{\tau,n}\end{math}, we define the 
weight of \begin{math}O\end{math}, \begin{math}W(O)\end{math}, to be the product of \begin{math}y\end{math} labels and 
the sign of \begin{math}O\end{math}, \begin{math}sgn(O)\end{math}, to be \begin{math}(-1)^{\ell(\mu)}\end{math}.
  For example, if $\tau = 13245$, then 
such an object \begin{math}O\end{math} constructed from the brick tabloid 
\begin{math}B= (2,8,3)\end{math} is pictured in Figure \ref{fig:On} where 
\begin{math}W(O) = y^7\end{math} and \begin{math}sgn(O) =(-1)^3\end{math}. 
It follows that 
\begin{equation}\label{eq:basic2}
n!\theta_{\tau}(h_n) = \sum_{O \in \mathcal{O}_{\tau,n}} sgn(O) W(O).
\end{equation}

Next we define a weight-preserving sign-reversing involution 
$I_{\tau}$ on \begin{math}\mathcal{O}_{\tau,n}\end{math}. Given an element 
\begin{math}O \in \mathcal{O}_{\tau,n}\end{math}, 
scan the cells of \begin{math}O\end{math} from left to right 
looking for the first 
cell \begin{math}c\end{math} such that either 

\begin{figure}[htbp]
  \begin{center}
    \includegraphics[width=0.6\textwidth]{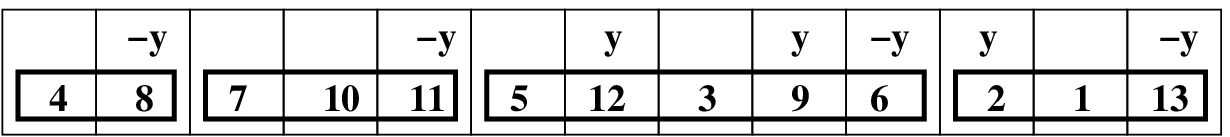}
    \caption{$I_{\tau}(O)$ for $O$ in Figure \ref{fig:On}.}
    \label{fig:IOn}
  \end{center}
\end{figure}

\begin{enumerate}
\item[(i)] $c$ is labeled with a $y$ or 

\item[(ii)] $c$ is a cell at the end of a brick \begin{math}b_i\end{math}, the number in cell \begin{math}c\end{math} 
is greater than the number in the first 
cell of the next brick \begin{math}b_{i+1}\end{math}, 
and there is no $\tau$-match 
in the cells of bricks \begin{math}b_i\end{math} and \begin{math}b_{i+1}\end{math}. 
\end{enumerate} 
In case (i), if \begin{math}c\end{math} is a cell in brick \begin{math}b_j\end{math}, then we split 
\begin{math}b_j\end{math} in to two bricks \begin{math}b_j^\prime\end{math} and \begin{math}b_j^{\prime \prime}\end{math} where 
\begin{math}b_j^\prime\end{math} contains all the cells of \begin{math}b_j\end{math} up to an including 
cell \begin{math}c\end{math} and \begin{math}b_j^{\prime \prime}\end{math} consists of the remaining cells 
of \begin{math}b_j\end{math} and we change the label on cell \begin{math}c\end{math} from \begin{math}y\end{math} to \begin{math}-y\end{math}. 
In case (ii), we combine the two bricks \begin{math}b_i\end{math} and \begin{math}b_{i+1}\end{math} into 
a single brick \begin{math}b\end{math} and change the label on cell \begin{math}c\end{math} from \begin{math}-y\end{math} to \begin{math}y\end{math}.  
For example, consider the element 
\begin{math}O \in \mathcal{O}_{13245,13}\end{math} pictured in Figure \ref{fig:On}.  Note that even 
though the number in the last cell of the first brick is greater than the the number in the first cell of the second brick, we cannot combine these two bricks because 
the numbers  \begin{math}4~8~7~10~11\end{math} would be a 13245-match. Thus the first 
place that we can apply the involution is on cell 5 which is labeled 
with a \begin{math}y\end{math} so  that \begin{math}I_{\tau}(O)\end{math} is the object pictured in Figure \ref{fig:IOn}. 

We claim that \begin{math}I_{\tau}\end{math} is an involution so that \begin{math}I_{\tau}^2\end{math} is the identity. 
To see this, consider case (i) where we split a brick \begin{math}b_j\end{math} at 
cell \begin{math}c\end{math} which is labeled with a \begin{math}y\end{math}.  
In that case, we let 
\begin{math}a\end{math} be the number in cell \begin{math}c\end{math} and \begin{math}a^\prime\end{math} be the number
in cell \begin{math}c+1\end{math} which must also be in brick \begin{math}b_{j}\end{math}. It must be the case 
that there is no cell labeled \begin{math}y\end{math} before cell \begin{math}c\end{math} since otherwise 
we would not use cell \begin{math}c\end{math} to define the involution. 
However, we  have to consider the possibility 
that when we spilt \begin{math}b_j\end{math} into \begin{math}b_j^\prime\end{math} and \begin{math}b_j^{\prime \prime}\end{math} 
that we might then be able to combine the brick \begin{math}b_{j-1}\end{math} with  
\begin{math}b_j^\prime\end{math} because the number in that last cell of 
\begin{math}b_{j-1}\end{math} is greater than the number in the first cell of 
\begin{math}b_j^\prime\end{math} and there is no $\tau$-match in the cells of 
\begin{math}b_{j-1}\end{math} and  \begin{math}b_j^\prime\end{math}.  
Since we always take an action on the left most cell possible 
when defining $I_\tau(O)$, we know that we cannot combine  \begin{math}b_{j-1}\end{math} and \begin{math}b_j\end{math} so that 
there must be a $\tau$-match in the cells of \begin{math}b_{j-1}\end{math} and \begin{math}b_j\end{math}. Clearly, 
that match must have involved the number \begin{math}a^\prime\end{math} and 
the number in cell $d$ which is the last cell in brick $b_{j-1}$. But that is impossible because then there would be two descents among the numbers
between cell $d$ and cell $c+1$ which would violate our assumption 
that $\tau$ has only one descent. Thus whenever we apply case (i) to 
define $I_\tau(O)$, 
the first action that we can take is to combine bricks 
\begin{math}b_j^\prime\end{math} and \begin{math}b_j^{\prime \prime}\end{math} so that \begin{math}I_{\tau}^2(O) = O\end{math}. 

If we are in case (ii), then again we can assume that there 
are no cells labeled \begin{math}y\end{math} that occur before cell \begin{math}c\end{math}.  When we combine 
brick \begin{math}b_i\end{math} and \begin{math}b_{i+1}\end{math}, then we will label cell \begin{math}c\end{math} with a \begin{math}y\end{math}. 
It is clear that combining the cells of \begin{math}b_i\end{math} and \begin{math}b_{i+1}\end{math} cannot 
help us combine the resulting brick \begin{math}b\end{math} with an earlier brick since 
it will be harder to have no $\tau$-matches with the larger brick \begin{math}b\end{math}. 
Thus the first place cell \begin{math}c\end{math} where we can apply the involution 
will again be cell \begin{math}c\end{math} which is now labeled with a \begin{math}y\end{math} so that \begin{math}I_{\tau}^2(O) =O\end{math} if we are in case (ii).

It is clear from our definitions that if \begin{math}I_{\tau}(O) \neq O\end{math}, 
then 
\begin{math}sgn(O)W(O) = - sgn(I_{\tau}(O))W(I_{\tau}(O))\end{math}. Hence it follows from 
(\ref{eq:basic2}) that 
\begin{equation}\label{eq:basic3} 
n!\theta_\tau(h_n) =  \sum_{O \in \mathcal{O}_{\tau,n}} sgn(O) W(O) = 
\sum_{O \in \mathcal{O}_{\tau,n}, I_{\tau}(O) =O} sgn(O) W(O). 
\end{equation}
Thus we must examine the fixed points of \begin{math}I_{\tau}\end{math}. So assume 
that \begin{math}O\end{math} is a fixed point of \begin{math}I_{\tau}\end{math}. 
First it is easy to see that there can be no cells which are labeled with 
\begin{math}y\end{math} so that numbers in each brick of \begin{math}O\end{math} must be increasing. 
Second we cannot combine two consecutive bricks \begin{math}b_i\end{math} and \begin{math}b_{i+1}\end{math} in 
\begin{math}O\end{math} which means either that there is an increase between the bricks 
\begin{math}b_i\end{math} and \begin{math}b_{i+1}\end{math} or there is a decrease between the bricks 
\begin{math}b_i\end{math} and \begin{math}b_{i+1}\end{math}, but there is a 
$\tau$-match in the cells of the bricks 
\begin{math}b_i\end{math} and \begin{math}b_{i+1}\end{math}. We claim that, in addition, the numbers in 
the first cells of the bricks must form an increasing sequence, 
reading from left to right. That is, suppose that 
\begin{math}b_i\end{math} and \begin{math}b_{i+1}\end{math} are two consecutive bricks in a fixed point \begin{math}O\end{math} of 
\begin{math}I_{\tau}\end{math} and that \begin{math}a > a^\prime\end{math} where \begin{math}a\end{math} is the number in the first 
cell of \begin{math}b_i\end{math} and \begin{math}a^\prime\end{math} is the number in the first cell of \begin{math}b_{i+1}\end{math}. 
Then clearly the number in the last cell of \begin{math}b_i\end{math} must be greater
than \begin{math}a^\prime\end{math} so that it must be the case that 
there is a $\tau$-match in the cells of \begin{math}b_i\end{math} and \begin{math}b_{i+1}\end{math}.  However 
\begin{math}a^\prime\end{math} is the least number that resides in the cells of \begin{math}b_i\end{math} and \begin{math}b_{i+1}\end{math} 
which means that the only way that \begin{math}a^\prime\end{math} could be part of 
a $\tau$-match that occurs in the cells of \begin{math}b_i\end{math} and \begin{math}b_{i+1}\end{math} is to 
have \begin{math}a^\prime\end{math} play the role of 1.  But 
since we are assuming that $\tau$ starts with 1, this would mean that if 
$a^\prime$ is part of a $\tau$-match, then that $\tau$-match must be  
entirely contained in \begin{math}b_{i+1}\end{math} which is impossible. Thus \begin{math}a^\prime\end{math} cannot 
be part of any $\tau$-match that occurs in the cells of \begin{math}b_i\end{math} and \begin{math}b_{i+1}\end{math}.
 However, this would mean that the $\tau$-match that occurs in the cells of \begin{math}b_i\end{math} and \begin{math}b_{i+1}\end{math} must either  be contained entirely in the cells of \begin{math}b_i\end{math} or entirely 
in the cells of \begin{math}b_{i+1}\end{math} which again is impossible.  Hence  it must be the case  that \begin{math}a < a^\prime\end{math}. 

Thus we have proved the following. 

\begin{lemma} \label{lemma:keytau}
Suppose that $\tau \in S_j$, $\tau$ starts with 1, and 
$\des{\tau} =1$.  Let $\theta_\tau:\Lambda \rightarrow 
\mathbb{Q}(y)$ be the  ring homomorphism defined on 
$\Lambda$ where $\mathbb{Q}(y)$ is the set of rational functions in 
the variable $y$ over the rationals $\mathbb{Q}$,  
$\theta_{\tau}(e_0) =1$ and  
$\theta_{\tau}(e_n) = \frac{(-1)^n}{n!} NM_{\tau,n}(1,y)$ for $n \geq 1$.
 Then 
\begin{equation}\label{keytau}
n!\theta_\tau (h_n) = \sum_{O \in \mathcal{O}_{\tau,n},I_\tau(O) = O}sgn(O)W(O)
\end{equation}
where $\mathcal{O}_{\tau,n}$ is the set of objects and 
$I_\tau$ is the involution defined above. Moreover, every 
fixed point $O$ of $I_\tau$ has the following three properties. 
\begin{enumerate}
\item There are no cells labeled with $y$ in $O$ so that the 
numbers in each brick of $O$ are increasing, 
\item the first numbers in each brick of $O$ form an increasing 
sequence, reading from left to right, and 
\item if \begin{math}b_i\end{math} and \begin{math}b_{i+1}\end{math} are two consecutive bricks in \begin{math}O\end{math}, then 
either (a) there is increase between \begin{math}b_i\end{math} and \begin{math}b_{i+1}\end{math} or (b) 
there is a decrease  between \begin{math}b_i\end{math} and \begin{math}b_{i+1}\end{math}, but there is 
$\tau$-match in the cells of \begin{math}b_i\end{math} and \begin{math}b_{i+1}\end{math}. 
\end{enumerate}
\end{lemma}

Now we specialize to the case of $\tau = 1324 \ldots p$ where 
$p \geq 5$. In this case,  we can make a finer analysis of the fixed 
points of $I_\tau$.  Let $O$ be a fixed point of $I_\tau$. 
By Lemma \ref{lemma:keytau}, we know that 1 is in the first cell of $O$. 
We claim that 2 must be in the second or third cell of $O$.  That is, 
suppose that $2$ is in cell $c$ where $c > 3$.  Then since 
there are no descents within any brick, $2$ must be the first cell 
of a brick. Moreover, since the minimal numbers in the bricks of $O$ form 
an increasing sequence, reading from left to right, $2$ must be in the first cell of
the second brick. Thus if $b_1$ and $b_2$ are the first two 
bricks in $O$, then 1 is in the first cell of $b_1$ and 2 is in the first cell of 
$b_2$.  But then we claim that there is no $\tau$-match in 
the cells of $b_1$ and $b_2$. That is, since $c > 3$, $b_1$ has 
at least three cells so that $O$ starts with an increasing sequence 
of length 3. 
But this means that $1$ cannot be part of a $\tau$-match. 
Similarly, no other cell of $b_1$ can be part of 
$\tau$-match because the 2 in cell $c$ is less than 
any of the remaining numbers of $b_1$. Thus if there is a 
$\tau$-match among the cells of $b_1$ and $b_2$, it would have 
to be entirely contained in $b_2$ which is impossible. But this 
would mean that we could apply case (ii) of the definition of 
$I_\tau$ to $b_1$ and $b_2$ which would violate our assumption 
that $O$ is a fixed point of $I_\tau$. Thus, we have two cases. \\
\ \\
{\bf Case 1.} 2 is in cell 2 of $O$.\\
\ \\
In this case there are two possibilities, namely, either 
(i) 1 and 2 are both in the first brick $b_1$ of $O$ or (ii) 
brick $b_1$ is a single cell filled with 1 and 2 is in the first cell of the second brick $b_2$ of 
$O$.  In either case, it is easy to see that 1 is not part of a 
$\tau$-match in $O$ and if we remove cell 1 from $O$ and 
subtract 1 from the numbers in the remaining cells, we would 
end up with a fixed point $O'$ of $I_\tau$ in $\mathcal{O}_{\tau,n-1}$. 
Now in case (i), it is easy to see that $sgn(O)W(O) = sgn(O')W(O')$ 
and in case (ii) since $b_1$ will have a label $-y$ on the first cell, 
$sgn(O)W(O) = (-y)sgn(O')W(O')$.  It follows that fixed points 
in Case 1 will contribute 
$(1-y)U_{\tau,n-1}(y)$ to $U_{\tau,n}(y)$.\\
\ \\
{\bf Case 2.} 2 is in cell 3 of $O$. \\
\ \\
Let $O(i)$ denote the number in cell  $i$ of $O$ and 
$b_1,b_2, \ldots $ be the bricks of $O$, reading from 
left to right. Since there are no descents within bricks in $O$
and the first numbers of each brick are increasing, reading 
from left to right, it must be the case that $2$ 
is in the first cell of $b_2$.  Thus $b_1$ has two cells. Note that 
$b_2$ must have 
 at least $p-2$ cells since otherwise, 
there could be no $\tau$-match contained in the cells of $b_1$ and 
$b_2$ and we could combine bricks $b_1$ and $b_2$ which would mean 
that $O$ is not a fixed point of $I_\tau$.  
But then the only reason that we cannot 
combine bricks $b_1$ and $b_2$ is that there is a $\tau$-match in 
the cells of $b_1$ and $b_2$ which could only start at position 1.

Next we claim 
that $O(p-1) =p-1$.  That is, since there is a $\tau$-match 
starting at position 1 and $p \geq 5$, we know that all 
the numbers in the first $p-2$ cells of $O$ are strictly 
less than $O(p-1)$. Thus $O(p-1) \geq p-1$.  Now if 
$O(p-1) > p-1$, then  let $i$ be least number
in the set $\{1,\ldots,p-1\}$ that is not contained in bricks $b_1$ and $b_2$.  Since 
the numbers in each brick are increasing and the minimal 
numbers of the bricks are increasing, 
the only possible position for $i$ is the first cell of brick $b_3$.  But then it follows that there is 
a decrease between bricks $b_2$ and $b_3$.  Since $O$ is a fixed 
point of $I_\tau$, this must mean that there is a $\tau$-match 
in the cells of $b_2$ and $b_3$. But since $\tau$ has only one 
descent, this $\tau$-match can only start at the cell $c$ which 
is the second to the last cell of $b_2$. Thus $c$ could be 
$p-1$ if $b_2$ has $p-2$ cells or $c > p-1$ if $b_2$ has
more than $p-2$ cells. In either case, 
$p-1 < O(p-1) \leq  O(c) < O(c+1) > O(c+2) =i$. But this is impossible 
since to have a $\tau$-match starting at cell $c$, we must 
have $O(c) < O(c+2)$.  Thus it must be the case that 
$O(p-1) = p-1$ and $\{O(1), \ldots ,O(p-1)\} = \{1, \ldots, p-1\}$. 

We now have two subcases.\\
\ \\
{\bf Case 2.a.} There is no $\tau$-match in $O$ starting at cell $p-1$.\\
\ \\
Then we claim that $O(p) =p$.  That is, if $O(p) \neq p$, then 
$p$ cannot be in $b_2$  so that $p$ must be in the first cell of the brick $b_3$. But then 
we claim that we could combine bricks $b_2$ and $b_3$.  That is, 
there will be a decrease between bricks $b_2$ and $b_3$ since 
$p < O(p)$ and $O(p)$ is in $b_2$. Since there is no $\tau$-match in 
$O$ starting at cell $p-1$, the only possible $\tau$-match among the cells 
of $b_2$ and $b_3$ would have to start at a cell $c \neq p-1$. But it can't be that  $c < p-1$ since then it would be the case that 
$O(c) < O(c+1) < O(c+2)$. Similarly, 
it cannot be that $c > p-1$ since then $O(c) > p$ and $p$ has 
 to be part of the $\tau$-match which is impossible 
since $O(c)$ must play 
 the role of 1 in the $\tau$-match. Thus it must be the case 
that $O(p) =p$. It then follows that 
if we let $O'$ be the result of removing the first $p-1$ cells from 
$O$ and subtracting $p-1$ from the remaining numbers, then $O'$ will be a fixed point of $I_\tau$ in $\mathcal{O}_{\tau,n-(p-1)}$.  Note that if $b_2$ has $p-2$ cells, then $O'$ will start with 
a brick with one cell and if $b_2$ has more than $p-2$ cells, then $O'$ will 
start with a brick with at least two cells. Since there is $-y$ coming 
from the brick $b_1$, it is easy to see that 
the fixed points in Case 2.a will contribute 
$-yU_{\tau,n-(p-1)}(y)$ to $U_{\tau,n}(y)$.\\
\ \\
{\bf Case 2.b.} There is a $\tau$-match starting cell $p-1$ in $O$.\\
\ \\
In this case, it must be that $O(p-1) < O(p) > O(p+1)$ so 
that $b_2$ must have $p-2$ cells and brick $b_3$ starts 
at cell $p+1$.  We claim that $b_3$ must have at least $p-2$ cells. 
That is, if $b_3$ has less than $p-2$ cells, then there could be no 
$\tau$-match among the cells of $b_2$ and $b_3$ so then 
we could combine $b_2$ and $b_3$ violating the fact that 
$O$ is a fixed point of $I_\tau$. 

In the general case, assume 
that in $O$, the bricks $b_2, \ldots, b_{k-1}$ all have $(p-2)$ cells. 
Then let  $r_1 =1$ and for $j=2,\ldots, k-1$, let 
$r_j=1+(j-1)(p-2)$. Thus $r_j$ is the position of 
the second to last cell of brick $b_j$ for $1 \le j \le k-1$. 
Furthermore, assume that there is a $\tau$-match starting at cell 
$r_j$ for $1\le j \le k-1$. It follows that 
$O(r_{k-1}) < O(r_{k-1}+1) > O(r_{k-1}+2)$ so that brick 
$b_k$ must start at cell $r_{k-1}+2$ and there is a decrease 
between bricks $b_{k-1}$ and $b_k$. But then it must be the case 
that $b_k$ has at least $p-2$ cells since if 
 $b_k$ has less than $p-2$ cells, we could combine bricks $b_{k-1}$ and 
$b_k$ violating the fact that $O$ is a fixed point of $I_\tau$.  
Let $r_k = 1+(k-1)(p-2)$ and assume that $O$ does not 
have a $\tau$-match starting at position $r_k$. 
Thus we have the situation pictured below.

$$
\underbrace{\begin{tabular}{|c|c|}
\multicolumn{1}{c}{}&
\multicolumn{1}{c}{}\\
\hline
$\scriptstyle{O(1)}$&$\scriptstyle{O(2)}$\\
\hline
\end{tabular}}_{b_1}
\underbrace{\begin{tabular}{|c|c|c|c|}
\multicolumn{1}{c}{}&
\multicolumn{1}{c}{}&
\multicolumn{1}{c}{$\scriptstyle{r_2}$}&
\multicolumn{1}{c}{}\\
\hline
$\scriptstyle{O(3)}$&$\scriptstyle{\cdots}$&$\scriptstyle{O(r_2)}$&$\scriptstyle{O(r_2+1)}$\\
\hline
\end{tabular}}_{b_2}
\scriptstyle{\ldots}
\underbrace{\begin{tabular}{|c|c|c|c|}
\multicolumn{1}{c}{}&
\multicolumn{1}{c}{}&
\multicolumn{1}{c}{$\scriptstyle{r_{k-1}}$}&
\multicolumn{1}{c}{}\\
\hline
~&$\scriptstyle{\cdots}$&$\scriptstyle{O(r_{k-1})}$&$\scriptstyle{O(r_{k-1}+1)}$\\
\hline
\end{tabular}}_{b_{k-1}}
\underbrace{\begin{tabular}{|c|c|c|c|c}
\multicolumn{1}{c}{}&
\multicolumn{1}{c}{}&
\multicolumn{1}{c}{$\scriptstyle{r_k}$}&
\multicolumn{1}{c}{}\\
\hline
~&$\scriptstyle{\cdots}$&$\scriptstyle{O(r_k)}$&$\scriptstyle{O(r_k+1)}$&$\scriptstyle{\cdots}$\\
\hline
\end{tabular}}_{b_k}
$$

First we claim $O(r_j) =r_j$  and 
$\{1,\ldots,r_j\}=\{O(1),\ldots,O(r_j)\}$ for $j=1, \ldots, k$. 
We have shown that $O(1) =1$ and that $O(r_2) = O(p-1) =p-1$ and 
$\{O(1), \ldots, O(p-1)\} = \{1, \ldots, p-1\}$. Thus assume by induction, 
$O(r_{j-1})=r_{j-1}$ and $\{1,\dots,r_{j-1}\}=\{O(1),\dots,O(r_{j-1})\}$. 
Since there is a $\tau$-match that starts at cell $r_{j-1}$ and $p \geq 5$, 
we know that all the numbers 
$$O(r_{j-1}), O(r_{j-1}+1), \ldots, O(r_{j-1}+p-3)$$
 are less 
than $O(r_j) = O(r_{j-1}+p-2)$. Since 
$\{1,\dots,r_{j-1}\}=\{O(1),\dots,O(r_{j-1})\}$, it follows that 
$O(r_j) \geq r_j$.  Next suppose that 
$O(r_j) > r_j$. Then let $i$ be the least number 
that does not lie in the bricks $b_1, \ldots, b_j$. Because 
the numbers in each brick increase and the minimal numbers 
in the bricks are increasing, it must be the case that 
$i$ is in the first cell of the next brick $b_{j+1}$. Now it cannot be 
that $j < k$ because then we have that 
$i = O(r_j+2) \leq r_j < O(r_j) < O(r_{j+1})$ which would violate 
the fact that there is a $\tau$-match in $O$ starting at cell $r_j$. 
If $j=k$, then it follows that there is 
a decrease between bricks $b_k$ and $b_{k+1}$ since $b_{k+1}$ starts with 
$i\leq r_k < O(r_k)$.  Since $O$ is a fixed 
point of $I_\tau$, this must mean that there is a $\tau$-match 
in the cells of $b_k$ and $b_{k+1}$. But since $\tau$ has only one 
descent, this $\tau$-match can only start at the cell $c$ which 
is the second to the last cell of $b_k$. Thus $c$ must be greater than 
$r_k$ because by hypothesis there cannot be a $\tau$-match starting at cell $r_k$. So $b_{k+1}$ must have more than
$p-2$ cells. In this case, we have that 
$i \leq r_k < O(r_k) \leq  O(c) < O(c+1) > O(c+2) =i$. But this cannot 
be since to have a $\tau$-match starting at cell $c$, we must 
have $O(c) < O(c+2)$.  Thus it must be the case that 
$O(r_j) = r_j$. But then it must be the case that 
$r_{j-1} =O(r_{j-1}) < O(d) < O(r_j) =r_j$ for $r_{j-1} < d < r_j$ so 
that $\{O(1), \ldots , O(r_j)\} = \{1,\ldots ,r_j\}$ as desired. 
Thus we have proved by induction that $O(r_j) =r_j$  and 
$\{1,\ldots,r_j\}=\{O(1),\ldots,O(r_j)\}$ for $j=1, \ldots, k$. 

This means that the sequence $O(1), \ldots, O(r_k)$ is completely 
determined. Next we claim that since there is no $\tau$-match 
starting at position $r_k$, it must be the case 
that $O(r_k+1) = r_k+1$. 
That is, if $O(r_k+1) \neq r_k+1$, then 
$r_k+1$ cannot be in brick $b_k$  so then $r_k+1$ must be in the first cell of the brick
$b_{k+1}$. But then 
we claim that we could combine bricks $b_k$ and $b_{k+1}$.  That is, 
there will be a decrease between bricks $b_k$ and $b_{k+1}$ since 
$r_k+1 < O(r_k+1)$ and $O(r_k+1)$ is in $b_k$. Since there is no $\tau$-match 
starting in $O$ at cell $r_k$, 
the only possible $\tau$-match among the cells 
of $b_k$ and $b_{k+1}$ would have to start at a cell $c \neq r_k$. 
Now  it cannot 
be that  $c < r_k$ since then $O(c) < O(c+1) < O(c+2)$. But  
it cannot be that $c > r_k$ since then $O(c) > r_k+1$ and $r_k+1$ would 
have to be part of the $\tau$-match which means that  $O(c)$ could 
not play the role of 1 in the $\tau$-match. Thus it must be the case 
that $O(r_k+1) =r_k+1$. It then follows that 
if we let $O'$ be the result of removing the first $r_k$ cells from 
$O$ and subtracting $r_k$ from each number in the 
remaining cells, then $O'$ will be a fixed point $I_\tau$ in $\mathcal{O}_{\tau,n-r_k}$.  Note that if $b_k$ has $p-2$ cells, then the first brick of $O'$ 
will have one cell and if $b_k$ has more than $p-2$ cells, then the first brick of $O'$ will have at least two cells. Since there is a factor $-y$ coming 
from each of the bricks $b_1, \ldots,b_{k-1}$, it is easy to see that 
the fixed points in Case 2.b will contribute 
$\sum_{k\geq 3} (-y)^{k-1} U_{\tau,n-((k-1)(p-2)+1)}(y)$ to $U_{\tau,n}(y)$.

Thus we have proved the following theorem.
\begin{theorem} \label{thm:132p}
Let $\tau=1324\dots p$ where 
$p \geq 5$. Then $\displaystyle NM_\tau(t,x,y)=\left(\frac1{U_\tau(t,y)}\right)^x$ where 
$\displaystyle U_\tau(t,y)=1+\sum_{n\geq1}U_{\tau,n}(y)\frac{t^n}{n!}$, 
$U_{\tau,1}(y)=-y$, and for $n \geq 2$,
$$U_{\tau,n}(y)=(1-y)U_{\tau,n-1}(y)+\sum_{k=1}^{\lfloor\frac{n-2}{p-2}\rfloor}(-y)^{}U_{\tau,n-(k(p-2)+1)}(y).$$
\end{theorem}

For example, we have computed the following.
\newpage
\noindent
$U_{13245,1}(y) = -y$,\\
$U_{13245,2}(y) = -y+y^2$,\\
$U_{13245,3}(y) = -y+2 y^2-y^3$,\\
$U_{13245,4}(y) = -y+3 y^2-3 y^3+y^4$,\\
$U_{13245,5}(y) = -y+5 y^2-6 y^3+4 y^4-y^5$,\\
$U_{13245,6}(y) = -y+7 y^2-12 y^3+10 y^4-5 y^5+y^6$,\\
$U_{13245,7}(y) = -y+9 y^2-21 y^3+23 y^4-15 y^5+6 y^6-y^7$,\\
$U_{13245,8}(y) = -y+11 y^2-34 y^3+47 y^4-39 y^5+21 y^6-7 y^7+y^8$,\\
$U_{13245,9}(y) = -y+13 y^2-51 y^3+88 y^4-90 y^5+61 y^6-28 y^7+8 y^8-y^9$,\\
$U_{13245,10}(y) = -y+15 y^2-72 y^3+153 y^4-189 y^5+156 y^6-90 y^7+36 y^8-9 y^9+y^{10}$,\\
$U_{13245,11}(y) = -y+17 y^2-97 y^3+250 y^4-368 y^5+361 y^6-252 y^7+127 y^8-45 y^9+10 y^{10}-y^{11}$.\\
\ \\
$U_{132456,1}(y) = -y$,\\
$U_{132456,2}(y) = -y+y^2$,\\
$U_{132456,3}(y) = -y+2 y^2-y^3$,\\
$U_{132456,4}(y) = -y+3 y^2-3 y^3+y^4$,\\
$U_{132456,5}(y) = -y+4 y^2-6 y^3+4 y^4-y^5$,\\
$U_{132456,6}(y) = -y+6 y^2-10 y^3+10 y^4-5 y^5+y^6$,\\
$U_{132456,7}(y) = -y+8 y^2-17 y^3+20 y^4-15 y^5+6 y^6-y^7$,\\
$U_{132456,8}(y) = -y+10 y^2-27 y^3+38 y^4-35 y^5+21 y^6-7 y^7+y^8$,\\
$U_{132456,9}(y) = -y+12 y^2-40 y^3+68 y^4-74 y^5+56 y^6-28 y^7+8 y^8-y^9$,\\
$U_{132456,10}(y) = -y+14 y^2-57 y^3+114 y^4-146 y^5+131 y^6-84 y^7+36 y^8-9 y^9+y^{10}$,\\
$U_{132456,11}(y) = -y+16 y^2-78 y^3+182 y^4-270 y^5+282 y^6-216 y^7+120 y^8-45 y^9+10 y^{10}-y^{11}$.\\
\ \\
$U_{1324567,1}(y) = -y$,\\
$U_{1324567,2}(y) = -y+y^2$,\\
$U_{1324567,3}(y) = -y+2 y^2-y^3$,\\
$U_{1324567,4}(y) = -y+3 y^2-3 y^3+y^4$,\\
$U_{1324567,5}(y) = -y+4 y^2-6 y^3+4 y^4-y^5$,\\
$U_{1324567,6}(y) = -y+5 y^2-10 y^3+10 y^4-5 y^5+y^6$,\\
$U_{1324567,7}(y) = -y+7 y^2-15 y^3+20 y^4-15 y^5+6 y^6-y^7$,\\
$U_{1324567,8}(y) = -y+9 y^2-23 y^3+35 y^4-35 y^5+21 y^6-7 y^7+y^8$,\\
$U_{1324567,9}(y) = -y+11 y^2-34 y^3+59 y^4-70 y^5+56 y^6-28 y^7+8 y^8-y^9$,\\
$U_{1324567,10}(y) = -y+13 y^2-48 y^3+96 y^4-130 y^5+126 y^6-84 y^7+36 y^8-9 y^9+y^{10}$,\\
$U_{1324567,11}(y) = -y+15 y^2-65 y^3+150 y^4-230 y^5+257 y^6-210 y^7+120 y^8-45 y^9+10 y^{10}-y^{11}$.\\
\ \\

Of course, one can use these initial values of the $U_{1324 \ldots p,n}(y)$ to 
compute the initial values of $NM_{1324 \ldots p}(t,x,y)$. 
For example, we have used Mathematica to compute the following initial terms 
of $NM_{13245}(t,x,y)$, $NM_{132456}(t,x,y)$, and $NM_{1324567}(t,x,y)$.

\begin{eqnarray*}
&&NM_{13245}(t,x,y) = 1+x y t+\frac{1}{2} \left(x y+x^2 y^2\right) t^2+
\frac{1}{6} \left(x y+x y^2+3 x^2 y^2+x^3 y^3\right) t^3+\\
&&\frac{1}{24} \left(x y+4 x y^2+7 x^2 y^2+x y^3+4 x^2 y^3+6 x^3 y^3+x^4 y^4\right) t^4+\\
&&\frac{1}{120}\left(x y+10 x y^2+15 x^2 y^2+11 x y^3+30 x^2 y^3+25 x^3 y^3+\right.\\
&& \left.x y^4+5 x^2 y^4+10 x^3 y^4+10 x^4 y^4+x^5 y^5\right) t^5+\\
&&\frac{1}{720} \left(x y+24 x y^2+31 x^2 y^2+62 x y^3+140 x^2 y^3+90 x^3 y^3+26 x y^4+91 x^2 y^4+\right.\\
&&\left.120 x^3 y^4+65 x^4 y^4+x y^5+6 x^2 y^5+15 x^3 y^5+20 x^4 y^5+15 x^5 y^5+x^6 y^6\right) t^6+\\
&&\frac{1}{5040}\left(x y+54 x y^2+63 x^2 y^2+273 x y^3+553 x^2 y^3+301 x^3 y^3+292 x y^4+\right.\\
&&840 x^2 y^4+875 x^3 y^4+350 x^4 y^4+57 x y^5+238 x^2 y^5+406 x^3 y^5+350 x^4 y^5+\\
&&\left.140 x^5 y^5+x y^6+7 x^2 y^6+21 x^3 y^6+35 x^4 y^6+35 x^5 y^6+21 x^6 y^6+x^7 y^7\right) t^7 + \\
&&\frac{1}{40320}\left(x y+116 x y^2+127 x^2 y^2+1068 x y^3+2000 x^2 y^3+966 x^3 y^3+2228 x y^4+\right. \\
&&5726 x^2y^4+5152 x^3 y^4+1701 x^4 y^4+1171 x y^5+4016 x^2 y^5+5474 x^3 y^5+3640 x^4 y^5+\\
&&1050 x^5 y^5+120 x y^6+575 x^2 y^6+1176 x^3 y^6+1316 x^4 y^6+840
x^5 y^6+266 x^6 y^6+\\
&&\left. x y^7+8 x^2 y^7+28 x^3 y^7+56 x^4 y^7+70 x^5 y^7+56 x^6 y^7+28 x^7 y^7+x^8 y^8\right) t^8+\cdots
\end{eqnarray*}

\begin{eqnarray*}
&&NM_{132456}(t,x,y) =1+x y t+\frac{1}{2} \left(x y+x^2 y^2\right) t^2+\frac{1}{6} \left(x y+x y^2+3 x^2 y^2+x^3 y^3\right) t^3+\\
&&\frac{1}{24} \left(x y+4 x y^2+7 x^2 y^2+x y^3+4 x^2 y^3+6 x^3 y^3+x^4 y^4\right) t^4+ \\
&&\frac{1}{120}\left(x y+11 x y^2+15 x^2 y^2+11 x y^3+30 x^2 y^3+25 x^3 y^3+\right.\\
&&\left. x y^4+5 x^2 y^4+10 x^3 y^4+10 x^4 y^4+x^5 y^5\right)t^5+ \\
&&\frac{1}{720} \left(x y+25 x y^2+31 x^2 y^2+66 x y^3+146 x^2 y^3+90 x^3 y^3+26 x y^4+91 x^2 y^4+\right.\\
&&\left.120 x^3 y^4+65 x^4 y^4+x y^5+6 x^2 y^5+15 x^3 y^5+20 x^4 y^5+15 x^5 y^5+x^6 y^6\right) t^6+\\
&&\frac{1}{5040}\left(x y+55 x y^2+63 x^2 y^2+297 x y^3+581 x^2 y^3+301 x^3 y^3+302 x y^4+\right.\\
&&868 x^2 y^4+896 x^3 y^4+350 x^4 y^4+57 x y^5+238 x^2 y^5+406 x^3 y^5+350 x^4 y^5+\\
&&\left.140 x^5 y^5+x y^6+7 x^2 y^6+21 x^3 y^6+35 x^4 y^6+35 x^5 y^6+21 x^6 y^6+x^7 y^7\right) t^7 + \\
&&\frac{1}{40320}\left(x y+117 x y^2+127 x^2 y^2+1153 x y^3+2092 x^2 y^3+966 x^3 y^3+2401 x y^4+\right.\\
&&6086 x^2
y^4+5348 x^3 y^4+1701 x^4 y^4+1191 x y^5+4096 x^2 y^5+5586 x^3 y^5+
3696 x^4 y^5+\\
&&1050 x^5 y^5+120 x y^6+575 x^2 y^6+1176 x^3 y^6+1316 x^4 y^6+840
x^5 y^6+266 x^6 y^6+\\
&&\left. x y^7+8 x^2 y^7+28 x^3 y^7+56 x^4 y^7+70 x^5 y^7+56 x^6 y^7+28 x^7 y^7+x^8 y^8\right) t^8+\cdots
\end{eqnarray*}

\begin{eqnarray*}
&&NM_{1324567}(t,x,y) = 1+x y t+\frac{1}{2} \left(x y+x^2 y^2\right) t^2+\frac{1}{6} \left(x y+x y^2+3 x^2 y^2+x^3 y^3\right) t^3+\\
&&\frac{1}{24} \left(x y+4 x y^2+7 x^2 y^2+x y^3+4 x^2 y^3+6 x^3 y^3+x^4 y^4\right) t^4+ \\
&&\frac{1}{120}\left(x y+11 x y^2+15 x^2 y^2+11 x y^3+30 x^2 y^3+25 x^3 y^3+\right.\\
&&\left. x y^4+5 x^2 y^4+10 x^3 y^4+10 x^4 y^4+x^5 y^5\right) t^5+\\
&&\frac{1}{720} \left(x y+26 x y^2+31 x^2 y^2+66 x y^3+146 x^2 y^3+90 x^3 y^3+26 x y^4+91 x^2 y^4+\right.\\
&&\left.120 x^3 y^4+65 x^4 y^4+x y^5+6 x^2 y^5+15 x^3 y^5+20 x^4 y^5+15 x^5 y^5+x^6 y^6\right) t^6+\\
&&\frac{1}{5040}\left(x y+56 x y^2+63 x^2 y^2+302 x y^3+588 x^2 y^3+301 x^3 y^3+302 x y^4+\right.\\
&&868 x^2 y^4+896 x^3 y^4+350 x^4 y^4+57 x y^5+238 x^2 y^5+406 x^3 y^5+350 x^4 y^5+\\
&&\left.140 x^5 y^5+x y^6+7 x^2 y^6+21 x^3 y^6+35 x^4 y^6+35 x^5 y^6+21 x^6 y^6+x^7 y^7\right) t^7+ \\
&&\frac{1}{40320}\left(x y+118 x y^2+127 x^2 y^2+1185 x y^3+2128 x^2 y^3+966 x^3 y^3+2416 x y^4+6126 x^2 y^4+\right.\\
&&5376 x^3 y^4+1701 x^4 y^4+1191 x y^5+4096 x^2 y^5+5586 x^3 y^5+3696 x^4 y^5+1050 x^5 y^5+\\
&&120 x y^6+575 x^2 y^6+1176 x^3 y^6+1316 x^4 y^6+840
x^5 y^6+266 x^6 y^6+x y^7+8 x^2 y^7+\\
&&\left. 28 x^3 y^7+56 x^4 y^7+70 x^5 y^7+56 x^6 y^7+28 x^7 y^7+x^8 y^8\right) t^8+ \cdots
\end{eqnarray*}

We note that there are many terms in these expansions which are 
easily explained. For example, we claim that for any $p \geq 4$, 
the coefficient 
of $x^ky^k$ in $NM_{1324 \ldots p,n}(x,y)$ is always the Stirling number 
$S(n,k)$ which is the number of set partitions of $\{1, \ldots,n\}$ 
into $k$ parts. That is, a permutation $\sg \in S_n$ that 
contributes to the coefficient $x^k y^k$ in $NM_{1324 \ldots p,n}(x,y)$ 
must have $k$ left-to-right minima and $k-1$ descents. Since 
each left-to-right minima of $\sg$ which is not the first element 
is always the second element of descent pair, it follows 
that if $1 = i_1 < i_2 < i_3 < \cdots < i_k$ are the positions of 
the left to right minima, then $\sg$ must be increasing 
in each of the intervals $[1,i_2),[i_2,i_3), \ldots, [i_{k-1},i_k),[i_k,n]$.
It is then easy to see that 
$$\{\sg_1, \ldots, \sg_{i_2-1}\}, 
\{\sg_{i_2}, \ldots, \sg_{i_3-1}\}, \ldots, 
\{\sg_{i_{k-1}}, \ldots, \sg_{i_k-1}\},\{\sg_{i_k}, \ldots, \sg_{n}\}$$
is just a set partition of $\{1, \ldots, n\}$ ordered by 
decreasing minimal elements. Moreover, it is easy to see 
that no such permutation can have a $1324 \ldots p$-match for any 
$p \geq 4$. Vice versa, if 
$A_1, \ldots, A_k$ is a set partition of 
$\{1, \ldots, n\}$ such that $min(A_1) > \cdots > min(A_k)$, then 
the permutation $\sg = A_k\uparrow A_{k-1}\uparrow \ldots A_1\uparrow$ 
is a permutation with $k$ left-to-right minima and $k-1$ descents 
where for any set $A \subseteq \{1, \ldots, n\}$, $A\uparrow$ is the 
list of the element of $A$ in increasing order. 
It follows that for any $p \geq 4$, 
\begin{enumerate}
\item  $NM_{1324 \ldots p,n}(x,y)|_{xy} =S(n,1) =1$, 
\item $NM_{1324 \ldots p,n}(x,y)|_{x^2y^2} =S(n,2) =2^{n-1}-1$, 
\item $NM_{1324 \ldots p,n}(x,y)|_{x^ny^n} =S(n,n) =1$, and 
\item $NM_{1324 \ldots p,n}(x,y)|_{x^ny^n} =S(n,n-1) =\binom{n}{2}$.
\end{enumerate}

We claim that  
$$NM_{1324 \ldots p,n}(x,y)|_{xy^2}= 
\begin{cases} 2^{n-1}-n & \mbox{if $n < p$ and } \\
2^{n-1}-n -(n-(p-1)) = 2^{n-1}-2n+p-1 & \mbox{if $n \geq p$}.
\end{cases}$$
That is, suppose that $\sg \in S_n$ contributes to $NM_{1324 \ldots p,n}(x,y)|_{xy^2}$. Then $\sg$ must have 1 left-to-right minima and one descent. 
It follows that $\sg$ must start with $1$ and have one descent. 
Now if $A$ is any subset of $\{2, \ldots, n\}$ and $B = \{2, \ldots, n\}-A$, 
then we let $\sg_A$ be the permutation $\sg_A = 1~A \uparrow  B\uparrow $. The only choices of $A$ that do not 
give rise to a permutation with one descent are $\emptyset$ and 
$\{2, \ldots, i\}$ for $i =2, \ldots, n$. It follows that there 
$2^{n-1}-n$ permutations that start with 1 and have 1 descent. 
Next consider when such a $\sg_A$ could have a $1234 \ldots p$-match. 
If the $1234 \ldots p$-match starts at position $i$, 
then it must be the case that $\red{\sg_i\sg_{i+1}\sg_{i+2}\sg_{i+3}} =1324$. 
This means that the only descent is at position $i+1$ and all the 
elements $\sg_{j}$ for $j \geq i+3$ are greater than or equal to $\sg_{i+3}$. 
But this means that all the elements between 1 and $\sg_{i+2}$ must 
appear in increasing order in $\sg_2 \ldots \sg_{i-1}$. It follows 
that $\sg_A$ is of the form $1 \ldots (q-2) q (q+2) (q+1) (q+2) \ldots n$. 
There are no such permutations if $n \leq p-1$ and there are 
$n-(p-1)$ such permutations if $n \geq p$ as $q$ 
can range from 1 to $n-(p-1)$.

We end this section by considering the 
special case where $\tau = 1324$ where the analysis of the 
fixed points of $I_{1324}$ is a bit different. 
Let $O$ be a fixed point of $I_{1324}$. 
By Lemma \ref{lemma:keytau}, we know that 1 is in the first cell of $O$. 
Again, we  claim that 2 must be in the second or third cell of $O$.  That is, 
suppose that $2$ is in cell $c$ where $c > 3$.  Then since 
there are no descents within any brick, $2$ must be in the first cell
of a brick. Moreover, since the minimal numbers in the bricks of $O$ form 
an increasing sequence, reading from left to right, $2$ must be in the first cell of 
the second brick. Thus if $b_1$ and $b_2$ are the first two 
bricks in $O$, then 1 is in the first cell of $b_1$ and 2 is in the first cell of 
$b_2$.  But then we claim that there is no ${1324}$-match in 
the elements of $b_1$ and $b_2$. That is, since $c > 3$, $b_1$ has 
at least three cells so that $O$ starts with an increasing sequence 
of length 3. 
But this means that $1$ can not be part of a $1324$-match. 
Similarly, no other cell of $b_1$ can be part of 
$1324$-match because the 2 in cell $c$ is smaller than 
any of the remaining numbers of $b_1$.  But this 
would mean that we could apply case (ii) of the definition of 
$I_{1324}$ to $b_1$ and $b_2$ which would violate our assumption 
that $O$ is a fixed point of $I_{1324}$. Thus, we have two cases. \\
\ \\
{\bf Case I.} 2 is in cell 2 of $O$.\\
\ \\
In this case there are two possibilities, namely, either 
(i) 1 and 2 lie in the first brick $b_1$ of $O$ or (ii) 
brick $b_1$ has one cell and 2 is the first cell of the second brick $b_2$ of 
$O$.  In either case, it is easy to see that 1 is not part of a 
$1324$-match and if we remove cell 1 from $O$ and 
subtract 1 from the elements in the remaining cells, we would 
end up with a fixed point $O'$ of $I_{1324}$ in $\mathcal{O}_{{1324},n-1}$. 
Now in case (i), it is easy to see that $sgn(O)W(O) = sgn(O')W(O')$ 
and in case (ii) since $b_1$ will have a label $-y$ on the first cell, 
$sgn(O)W(O) = (-y)sgn(O')W(O')$.  It follows that fixed points 
in Case 1 will contribute 
$(1-y)U_{{1324},n-1}(y)$ to $U_{{1324},n}(y)$.\\
\ \\
{\bf Case II.} 2 in cell 3 of $O$. \\
\ \\
Let $O(i)$ denote the element in $i$ cell of $O$ and 
$b_1,b_2, \ldots $ be the bricks of $O$, reading from 
left to right. Since there are no descents within bricks in $O$ 
and the minimal elements in the bricks are increasing,  
we know that $2$ is in the first cell of a brick $b_2$.   
Thus $b_1$ has two cells. But then $b_2$ must have at least two cells since if $b_2$ has one cell,  
there could be no ${1324}$-match contained in the cells of $b_1$ and 
$b_2$ and we could combine bricks $b_1$ and $b_2$ which would mean 
that $O$ is not a fixed point of $I_{1324}$. Thus $b_1$ has two cells and  $b_2$ has at least two cells. 
But then the only reason that we could not 
combine bricks $b_1$ and $b_2$ is that there is a ${1324}$-match in 
the cells of $b_1$ and $b_2$ which could only start at the first cell.

We now have two subcases.\\
\ \\
{\bf Case II.a.} There is no ${1324}$-match in $O$ starting at cell $3$.\\
\ \\
Then we claim that $\{O(1),O(2),O(3),O(4)\}=\{1,2,3,4\}$.  
That is, if $\{O(1),O(2),O(3),O(4)\}\neq \{1,2,3,4\}$, then  
let $i = min(\{1,2,3,4\} -\{O(1),O(2),O(3),O(4)\})$. Since 
there is a 1324-match starting at position 1, it follows 
that $O(4) > 4$ since $O(4)$ is the fourth largest element 
in $\{O(1),O(2),O(3),O(4)\}$. Since the minimal elements 
of the bricks of $O$ are increasing, it must be that $i$ is 
the first element in brick $b_3$. But then 
we claim that we could combine bricks $b_2$ and $b_3$.  That is, 
there will be a decrease between bricks $b_2$ and $b_3$ since 
$i < O(4)$ and $O(4)$ is in $b_2$. Since there is no ${1324}$-match in 
$O$ starting at cell $3$, the only possible ${1324}$-match among the elements 
in $b_2$ and $b_3$ would have start at a cell $c > 3$. But then 
$O(c) > i$, which is impossible since it would have to play the 
role of 1 in the ${1324}$-match and $i$ would have to play the role of 2 in 
the ${1324}$-match since $i$ occupies the first cell of $b_3$. 
Thus it must be the case that $O(1) =1$, $O(2) =3$, $O(3) = 2$, and 
$O(4) = 4$. 

It then follows that 
if we let $O'$ be the result of removing the first $3$ cells from 
$O$ and subtracting $3$ from the remaining elements, then $O'$ will be a fixed point $I_{1324}$ in $\mathcal{O}_{{1324},n-3}$.  Since there is $-y$ coming 
from the brick $b_1$, it is easy to see that 
the fixed points in Case II.a will contribute 
$-yU_{{1324},n-3}(y)$ to $U_{{1324},n}(y)$.\\
\ \\
{\bf Case II.b.} There is a ${1324}$-match starting a $3$ in $O$.\\
\ \\
In this case, it must be that $O(3) < O(4) > O(5)$ so 
that $b_2$ must have two cells and brick $b_3$ starts 
at cell $5$.  We claim that $b_3$ must have at least two cells. 
That is, if $b_3$ has one cell, then there could be no 
${1324}$-match among the cells of $b_2$ and $b_3$ so that 
we could combine $b_2$ and $b_3$ violating the fact that 
$O$ is a fixed point of $I_{1324}$. 

In the general case, assume 
that in $O$, the bricks $b_2, \ldots, b_{k-1}$ all have two cells and 
there are $1324$-matches starting at cells $1,3, 
\ldots, 2k-3$ but there is no $1324$-match starting at cell
$2k-1$ in $O$.  Then we know that $b_k$ has least two cells. 
Let $c_i < d_i$ be the numbers in the first two cells of brick $b_i$ for 
$i =1, \ldots, k$. Then we have that $\red{c_id_ic_{i+1}d_{i+1}}=1324$ for $1\leq i \leq k-1$. This means that $c_i<c_{i+1}<d_i<d_{i+1}$. 

First we claim that it must be the case that 
$\{O(1),\ldots, O(2k)\} =\{1, \ldots, 2k\}$. If not 
there is a number greater than $2k$ that occupies one of the first $2k$ cells. Let $M$ be the greatest such number. If $M$ occupies one of the first $2k$ cells then there must be a number less than $2k$ that occupies one of the last $n-2k$ cells. Let $m$ be the least such number. Since numbers in bricks are increasing, $M$ must occupy the last cell in one of the first $k-1$ bricks or occupy cell $2k$. If $M$ occupies the last cell in one of the first $k-1$ bricks, then $M$ is part of a $\tau$-match 
$$\dots~
\begin{tabular}{|c|c|}
\hline
$c_i$&$M$\\
\hline
\end{tabular}~
\begin{tabular}{|c|c|}
\hline
$c_{i+1}$&$d_{i+1}$\\
\hline
\end{tabular}~\ldots .
$$
But then 
$\red{c_i~M~c_{i+1}~d_{i+1}}=1~3~2~4$ implies that $M<d_{i+1}$ which 
contradicts our choice of $M$ as the greatest number in the first $2k$ cells.
Thus $M$ cannot  occupy the last cell in one of the first $k-1$ bricks. This means that $M$ must occupy cell $2k$ in $O$.

Since numbers in bricks are increasing, $m$ must occupy the first cell of $b_{k+1}$. But then there 
is a descent between bricks $b_k$ and $b_{k+1}$ so 
that $m$ must be part of a $1324$-match.  But the only way 
this can happen is if in the $1324$-match involving $m$, 
$m$ plays the role of 
$2$ and the numbers in the last two cells of brick $b_k$ play the role of $1~3$.
Since, we are assuming that a $1324$-match does not start at cell $2k-1$ which is the cell that the number $c_k$ occupies, the numbers in the last two cells of 
brick $b_k$ must be greater than or equal to $d_k =M$ which is impossible 
since $m < M$. Thus it must be the case that 
$\{O(1),\ldots, O(2k)\} =\{1, \ldots, 2k\}$ and that $d_k =2k$. 
It now follows that if we remove the first $2k-1$ cells 
from $O$ and replace each remaining number $i$ in $O$ by 
$i-(2k-1)$, then we will end up with a fixed point in 
$O'$ of $I_{1324}$ in $\mathcal{O}_{n-(2k-1)}$. Thus each 
such fixed point $O$ will contribute 
$(-y)^{k-1}U_{n-2k+1}(y)$ to $U_n(y)$.  

\begin{figure}[htbp]
  \begin{center}
    \includegraphics[width=0.5\textwidth]{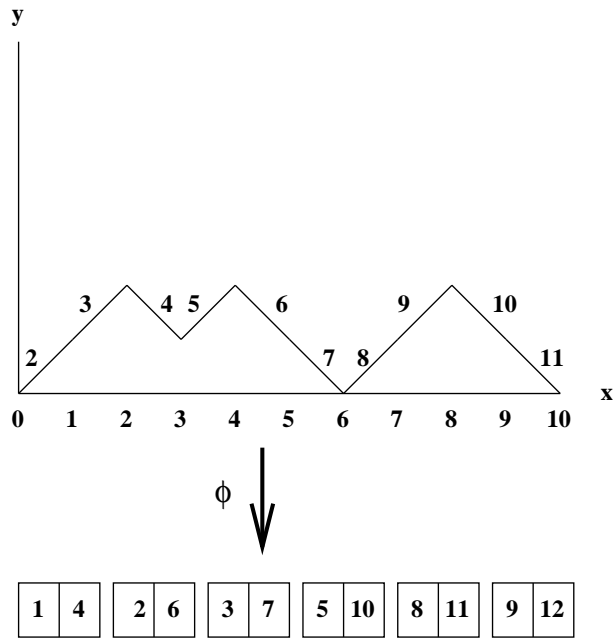}
    \caption{The bijection $\phi$.}
    \label{fig:DCbij}
  \end{center}
\end{figure}

The only thing left 
to do is to count the number of such fixed points $O$. That is, 
we must count the number of sequences 
$c_1 d_1 c_2 d_2 \ldots c_k d_k$ such that 
(i) $c_1 =1$, (ii) $c_2=2$, (iii) $d_k =2k$, 
(iv) $\{c_1,d_1, \ldots, c_k,d_k\} = \{1,2, \ldots, 2k\}$, and 
(v)  $\red{c_i d_i c_{i+1} d_{i+1}} = 1324$ for each $1 \leq i \leq k-1$. 
We claim that there are $C_{k-1}$ such sequences where 
$C_n =\frac{1}{n+1} \binom{2n}{n}$ is the $n$-th Catlan number. It is well 
known that $C_{k-1}$ counts the number of Dyck paths of length 
$2k-2$. A Dyck path of length $2k-2$ is a path that starts at $(0,0)$ and ends at $(2k-2,0)$ 
and consists of either {\em up-steps} (1,1) or {\em down-steps} (1,-1) in such a way that the 
path never goes below the $x$-axis. Thus we will give 
a bijection $\phi$ between the set of Dyck paths of length $2k-2$ and the 
set of sequences $c_1,d_1, \ldots, c_k,d_k$ satisfying conditions 
(i)-(v). The map $\phi$ is quite simple. That is, 
suppose that we start with a Dyck path $P = (p_1,p_2, \ldots, p_{2k-2})$ 
of length $2k-2$. First, label 
the segments $p_1, \ldots, p_{2k-2}$ with $2, \ldots, 2k-1$, respectively. 
Then $\phi(P)$ is the sequence $c_1 d_1 \ldots c_k d_k$ where 
$c_1 =1$ and $c_2 \ldots c_k$ are the labels of the up-steps of 
$P$, reading from left to right, $d_1 \ldots d_{k-1}$ are the 
labels of the down steps, reading from left to right, and 
$d_{2k}=2k$. We have pictured an example in Figure \ref{fig:DCbij}
of the bijection $\phi$ in the case where $k=6$. 

It is easy to see by construction that if $P$ is a Dyck path of 
length $2k-2$ and 
$\phi(P) = c_1 d_1 \ldots c_k d_k$, then 
$c_1 <  c_2 < \cdots < c_k$ and $d_1 <  d_2 < \cdots < d_k$. 
Moreover, since each Dyck path must start with an up-step, we have 
that $c_2 =2$. Clearly $c_1 =1$, $d_k=2k$, and 
$\{c_1,d_1, \ldots, c_k,d_k\} = \{1,\ldots, 2k\}$ by construction. 
Thus $c_1 d_1 \ldots c_k d_k$ satisfies conditions (i)-(iv). 
For condition (v), note that $c_1 = 1 < d_1 > 2 =c_2 <d_2$ so 
that $\red{c_1 d_1 c_2 d_2} = 1~3~2~4$. If $2 \leq i \leq k-1$, 
then note that $c_i$ equals the label of the $(i-1)$st up-step, 
$c_{i+1}$ equals the label of the $i$-th up-step, and $d_i$ is 
the label of $i$-th down-step. Since in a Dyck path, the 
$i$-th down-step must occur after the $i$-th up-step, it 
follows that $c_i < c_{i+1} < d_i < d_{i+1}$ so that 
$\red{c_i d_i c_{i+1} d_{i+1}} = 1~3~2~4$. Vice versa, 
if we start with a sequence $c_1 d_1 \ldots c_k d_k$ satisfying 
conditions (i)-(v) and create a path $P= (p_1, \ldots, p_{2k-2})$ with labels 
$2, \ldots, 2k-1$ such that $p_j$ is an up-step if 
$j+1 \in \{c_2, \ldots, c_{k}\}$ and $p_j$ is an down-step if 
$j+1 \in \{d_1, \ldots, d_{k-1}\}$, then condition (iii) ensures 
$P$ starts with an up-step and condition (v) ensures that 
the $i$-th up-step occurs before the $i$-th down step so that 
$P$ will be a Dyck path. Thus $\phi$ is a bijection between the set of Dyck paths of length $2k-2$ and the 
set of sequence $c_1,d_1, \ldots, c_k,d_k$ satisfying conditions 
(i)-(v).

It follows that fixed points $O$ of $I_{1324}$ where 
the bricks $b_1, b_2, \ldots, b_{k-1}$ are of size $2$ and 
there are $1324$-matches starting at positions $1,3, 
\ldots, 2k-3$ in $O$, but there is no $1324$-match starting at position 
$2k-1$ in $O$ contribute 
$C_{k-1} (-y)^{k-1} U_{n-2k+1}(y)$ to $U_n(y)$. Thus we have proved the 
following theorem. 

\begin{theorem}\label{thm:1324}
$\displaystyle NCM_{1324}(t,x,y)=\left(\frac1{U_{1324}(t,y)}\right)^x$ where $\displaystyle U_{1324}(t,y)=1+\sum_{n\geq1}U_{1324,n}(y)\frac{t^n}{n!}$, 
$U_{1324,1}(y)=-y$, and for $n\geq 2$,
$$U_{1324,n}(y)=(1-y)U_{1234,n-1}(y)+\sum_{k=2}^{\lfloor\frac{n}{2}\rfloor}(-y)^{k-1}C_{k-1}U_{1324,n-2k+1}(y)$$
where $C_k$ is the $k^{th}$ Catalan number.
\end{theorem}

In this case, one can easily compute that \\
\ \\
$U_{1324,1}(y) = -y$,\\
$U_{1324,2}(y) = -y+y^2$,\\
$U_{1324,3}(y) = -y+2 y^2-y^3$,\\
$U_{1324,4}(y) = -y+4 y^2-3 y^3+y^4$,\\
$U_{1324,5}(y) = -y+6 y^2-8 y^3+4 y^4-y^5$,\\
$U_{1324,6}(y) = -y+8 y^2-18 y^3+13 y^4-5 y^5+y^6$,\\
$U_{1324,7}(y) = -y+10 y^2-32 y^3+36 y^4-19y^5+6 y^6-y^7$,\\
$U_{1324,8}(y) = -y+12 y^2-50 y^3+85 y^4-61 y^5+26 y^6-7 y^7+y^8$,\\
$U_{1324,9}(y) = -y+14 y^2-72 y^3+166 y^4-170 y^5+94 y^6-34 y^7+8 y^8-y^9$,\\
$U_{1324,10}(y) = -y+16 y^2-98 y^3+287 y^4-412y^5+296 y^6-136 y^7+43 y^8-9 y^9+y^{10}$,\\
$U_{1324,11}(y) = -y+18 y^2-128 y^3+456 y^4-854 y^5+824 y^6-473 y^7+188 y^8-53 y^9+10 y^{10}-y^{11}$.\\

This, in turn, allows us to compute the first few terms 
of the generating function $NM_{1324}(t,x,y)$. That is, 
one can use Mathematica to compute that 
\begin{eqnarray*}
&&NM_{1324}(t,x,y) = \\
&&1+t x y+\frac{1}{2} t^2 \left(x y+x^2 y^2\right)+\frac{1}{6} t^3 \left(x y+x y^2+3 x^2 y^2+x^3 y^3\right)+\\
&&\frac{1}{24} t^4 \left(x y+3
x y^2+7 x^2 y^2+x y^3+4 x^2 y^3+6 x^3 y^3+x^4 y^4\right)+\\
&&\frac{1}{120} t^5 \left(x y+9 x y^2+15 x^2 y^2+8 x y^3+25 x^2 y^3+25 x^3 y^3+x y^4+ \right.\\
&&\left. 5 x^2
y^4+10 x^3 y^4+10 x^4 y^4+x^5 y^5\right)+\\
&&\frac{1}{720} t^6 \left(x y+23 x y^2+31 x^2 y^2+47 x y^3+119 x^2 y^3+90 x^3 y^3+20 x y^4+73 x^2 y^4+\right.\\
&&\left. 105
x^3 y^4+65 x^4 y^4+x y^5+6 x^2 y^5+15 x^3 y^5+20 x^4 y^5+15 x^5 y^5+x^6 y^6\right)+\\
&&\frac{1}{5040}t^7 \left(x y+53 x y^2+63 x^2 y^2+221 x y^3+490
x^2 y^3+301 x^3 y^3+202 x y^4+\right.\\
&&\left.637 x^2 y^4+749 x^3 y^4+350 x^4 y^4+47 x y^5+196 x^2 y^5+343 x^3 y^5+315 x^4 y^5+\right.\\
&&\left. 140 x^5 y^5+x y^6+7 x^2 y^6+21 x^3
y^6+35 x^4 y^6+35 x^5 y^6+21 x^6 y^6+x^7 y^7\right) + \\
&& \frac{1}{40320}\left(x y+115 x y^2+127 x^2 y^2+922 x y^3+1838 x^2 y^3+966 x^3 y^3+1571 x y^4+\right.\\
&&4421 x^2
y^4+4466 x^3 y^4+1701 x^4 y^4+795 x y^5+2890 x^2 y^5+4270 x^3 y^5+
3164 x^4 y^5+\\
&&1050 x^5 y^5+105 x y^6+495 x^2 y^6+1008 x^3 y^6+1148 x^4 y^6+770 x^5
y^6+266 x^6 y^6+x y^7+\\
&&\left. 8 x^2 y^7+28 x^3 y^7+56 x^4 y^7+70 x^5 y^7+56 x^6 y^7+28 x^7 y^7+x^8 y^8\right) t^8+
\cdots 
\end{eqnarray*}

We note that there are other methods to compute $NM_{1324}(t,1,1)$. 
That is, Elizalde \cite{E06} developed recursive techniques to 
find the coefficients of the series $NM_{1324}(t,1,1)$.

\section{Permutations with no $1324 \ldots p$-matches and one or two descents}

In this section, we will show how we can use 
Theorem \ref{thm:132p} and Theorem \ref{thm:1324} to find the 
generating function for the number of permutations 
$\sg \in S_n$ which have no $1324 \ldots p$-matches and have exactly 
$k$ descents for $k=1,2$ and 
$p \geq 4$.

That is, fix $p \geq 4$ and let 
$d^{(i)}_{n,p}$ denote the number of $\sigma \in S_n$ such 
that $1324\ldots p\mbox{-mch}(\sg) =0$ and $\des{\sg}=i$. 
Our goal is to compute 
$$D_p^{(i)}(t) = \sum_{n \geq 0}d^{(i)}_{n,p} \frac{t^n}{n!} = NM_{1324\ldots p}(t,1,y)|_{y^{i+1}}$$
for $i=1$ and $i=2$. 

To this end, we first want to compute 
$U_{1324 \ldots p,n}(y)|_{y}$, $U_{1324 \ldots p,n}(y)|_{y^2}$, and 
$U_{1324 \ldots p,n}(y)|_{y^3}$. That is, we want to 
compute the number of fixed points of $I_{1324 \ldots p}$ 
that have either 1, 2, or 3 bricks. 
Clearly there is only one fixed point of $I_{1324 \ldots p}$ of 
length $n$ which has just one brick since in that case, the 
underlying permutation must be the identity. In such a situation, 
the last cell of the brick is labeled with $-y$ so that 
for all $n \geq 1$ and all $p \geq 4$, $U_{1324 \ldots p,n}(y)|_{y} =-1$. 
Hence 
\begin{equation}\label{eq:Uygf}
U_{1324 \ldots p}(t,y)|_{y} = 1-e^t.
\end{equation}

Next we consider the fixed points of $I_{1324 \ldots p}$ which  are 
of length $n$ and consists 
of two bricks, a brick $B_1$ of length $b_1$ followed by  
a brick $B_2$ of length $b_2$. Note that in this case, 
the last cells of $B_1$ and $B_2$ are labeled with $-y$ so 
that the weight of all such fixed points is $y^2$. 
Suppose the underlying 
permutation is $\sg = \sg_1 \ldots \sg_n$.  Then there 
are two cases.\\
{\bf Case 1.}  There is an increase between the two bricks, i.e. 
$\sg_{b_1} < \sg_{b_1+1}$.\\

In this case, it easy to see that $\sg$ must be the identity permutation 
and, hence, there are $n-1$ fixed points in case 1 since 
$b_1$ can range from $1$ to $n-1$. \\
\ \\
{\bf Case 2.}  There is an decrease between the two bricks, i.e. 
$\sg_{b_1} > \sg_{b_1+1}$.\\

In this case, there must be a $1324 \ldots p$-match in 
the elements in the bricks of $B_1$ and $B_2$ which means 
that it must be the case that 
$\red{\sg_{b_1-1} \sg_{b_1} \sg_{b_1+1} \ldots \sg_{b_1+p-2}} = 1324 \ldots p$. Now suppose that $\sg_{b_1-1} =x$. Since $\sg_{b_1+1}$ is the smallest 
element in brick $B_2$ and the elements in brick $B_2$ increase and 
$\sg_{b_1} > \sg_{b_1+1}$, it must be the case that 
$1, \ldots,x-1$ must lie in brick $B_1$. It cannot be that 
$\sg_{b_1} = x+1$ since $\sg_{b_1} > \sg_{b_1+1}$. Thus it must be the 
case that 
$\sg_{b_1+1} =x+1$ and $\sg_{b_1} = x+2$ since $\sg_{b_1} < \sg_{b_1+2}$. 
Thus brick $B_1$ consists of the elements $1, \ldots, x, x+2$. Hence 
there are $n-p+1$ possibilities in case 2 if $n \geq p$ and no possibilities 
in case 2 if $n < p$. \\
\ \\

It follows that 
\begin{equation}
U_{1324 \ldots p,n}|_{y^2} = \begin{cases} 0 & \mbox{if } n =0,1 \\
n-1 & \mbox{if } 2 \leq n \leq p-1 \ \mbox{and} \\
2n-p & \mbox{if } n \geq p.
\end{cases}
\end{equation}
Note that 
\begin{eqnarray*}
\sum_{n \geq p} (2n-p) \frac{t^n}{n!} &=& 2t \sum_{n \geq p} \frac{t^{n-1}}{(n-1)!} - p \sum_{n \geq p} \frac{t^n}{n!} \\
&=& 2t(e^t - \sum_{n=0}^{p-2} \frac{t^n}{n!}) -p(e^t  -\sum_{n=0}^{p-1} \frac{t^n}{n!}) \\
&=& (2t-p)e^t + p + \sum_{n=1}^{p-1}(p-2n) \frac{t^n}{n!}.
\end{eqnarray*}
Thus 
\begin{eqnarray}\label{eq:Uy2gf}
U_{1324 \ldots p}(t,y)|_{y^2} &=&  (2t-p)e^t + p + \sum_{n=1}^{p-1} (p-2n) \frac{t^n}{n!} + \sum_{n=2}^{p-1} (n-1) \frac{t^n}{n!} \nonumber \\
&=& (2t-p)e^t + p + \sum_{n=1}^{p-1} (p-n-1) \frac{t^n}{n!}.
\end{eqnarray}

Next we consider a fixed point  of 
$I_{1324 \ldots p}$ which has 3 bricks, $B_1$ of size $b_1$ followed by  
$B_2$ of size $b_2$  followed by $B_3$ of size $b_3$. Let 
$\sg = \sg_1 \ldots \sg_n$ be the underlying  permutation. The weight of all 
such fixed points is $-y^3$. We then have 4 cases. \\
\ \\
{\bf Case a.} There are increases between $B_1$ and $B_2$ and 
between $B_2$ and $B_3$, i.e. $\sg_{b_1} < \sg_{b_1+1}$ and 
$\sg_{b_1+b_2} < \sg_{b_1+b_2+1}$. \\

In this case, it is easy to see that $\sg$ must be the identity permutation 
so that there are $\binom{n-1}{2}$ possibilities in case 1 if $n \geq 3$. \\
\ \\
{\bf Case b.} There is an increase between $B_1$ and $B_2$ and a decrease 
between $B_2$ and $B_3$, i.e. $\sg_{b_1} < \sg_{b_1+1}$ and 
$\sg_{b_1+b_2} > \sg_{b_1+b_2+1}$. \\
 
In this case, it must be the case that $\sg_1 < \cdots < \sg_{b_1+b_2}$ and 
$$\red{\sg_{b_1+b_2-1}\sg_{b_1+b_2}\sg_{b_1+b_2+1} 
\ldots \sg_{b_1+b_2 +p-2}} = 
1324 \ldots p.$$  Then we can argue exactly as in case 2 above 
that there must exist an $x$ such that $\sg_{b_1+b_2-1} =x$ and 
$1, \ldots, x-1$ must occur to the left of $\sg_{b_1+b_2-1}$, 
$\sg_{b_1+b_2} =x+2$ and $\sg_{b_1+b_2+1} =x+1$.  Then for any fixed 
$x\geq 2$, 
we have $x-1$ choices for the length of $B_1$ so that 
we have 
$\sum_{x=2}^{n-p+1}(x-1) = \binom{n-p+1}{2}$ possibilities if $n \geq p+1$ 
and no possibilities if $n \leq p$. \\
\ \\
{\bf Case c.} There is a decrease between $B_1$ and $B_2$ and an increase 
between $B_2$ and $B_3$, i.e. $\sg_{b_1} > \sg_{b_1+1}$ and 
$\sg_{b_1+b_2} < \sg_{b_1+b_2+1}$. \\
 
In this case, it must be the case that $\sg_{b_1+1}< \cdots < \sg_{n}$ and 
$\red{\sg_{b_1-1}\sg_{b_1}\sg_{b_1+1} \ldots \sg_{b_1+p-2}} = 
1324 \ldots p$.  Again we can argue as in case 2 above 
that there must exist an $x$ such that $\sg_{b_1-1} =x$ and 
$1, \ldots, x-1$ must occur to the left of $\sg_{b_1-1}$, 
$\sg_{b_1} =x+2$ and $\sg_{b_1+1} =x+1$.  Then for any fixed $x$, 
we have $n-1-(x+p-1)$ choices for the length of $B_2$ so that 
we have $\sum_{x=1}^{n-p} n-x-p-1 = \binom{n-p+1}{2}$ 
possibilities if $n \geq p+1$ and no possibilities if $n \leq p$. \\
\ \\
{\bf Case d.} There are decreases between $B_1$ and $B_2$ and 
between $B_2$ and $B_3$, i.e. $\sg_{b_1} > \sg_{b_1+1}$ and 
$\sg_{b_1+b_2} > \sg_{b_1+b_2+1}$. \\
 
In this case, there  are two subcases. \\
{\bf Subcase d.1} $p=4$.\\

Now we must have  
$\red{\sg_{b_1-1}\sg_{b_1}\sg_{b_1+1} \sg_{b_1+2}} = 
1324$.  First suppose that $b_2 =2$. Then we also have 
that $\red{\sg_{b_1+1}\sg_{b_1+2}\sg_{b_1+3} \sg_{b_1+4}} = 
1324$. It follows that $x= \sg_{b_1-1} < \sg_{b_1+1}< \sg_{b_1+3}$. 
Since $\sg_{b_1+1}$ is the smallest element in brick $B_2$ and 
$\sg_{b_1+3}$ is the smallest element in brick $B_3$, it must 
be the case that $1, \ldots x-1$ lie in brick $B_1$ and 
that $\sg_{b_1+1} =x+1$.  It also must be the case 
that $\sg_{b_1} < \sg_{b_1+2} < \sg_{b_1+4}$ so that 
$\sg_{b_1}, \sg_{b_1+2} \in \{x+2,x+3\}$ and $\sg_{b_1+4}=x+4$. 
Thus there are two possibilities for each $x$. As $x$ can vary between 
1 and $n-5$ in this case, we have $2(n-5)$ possibilities 
if $b_2=2$ and $n \geq 6$ and no possibilities if $n <6$.

Next consider the case where $b_2 \geq 3$. 
Again we must have  
$\red{\sg_{b_1-1}\sg_{b_1}\sg_{b_1+1} \sg_{b_1+2}} = 
1324$. Similarly we must have $\red{\sg_{b_1+b_2-1}\sg_{b_1+b_2}\sg_{b_1+b_2+1} \sg_{b_1+b_2+2}} = 
1324$, but this condition does not involve $\sg_{b_1+1}$. Nevertheless, 
these conditions still force that $\sg_{b_1} < \sg_{b_1+b_2+1}$ so 
that if $\sg_{b_1-1} =x$, then $x$ is less than the least elements 
in bricks $B_2$ and $B_3$ so that $1, \ldots, x-1$ must be in 
brick $B_1$ and $\sg_{b_1+1} = x+1$.  However in this case, 
$\red{\sg_{b_1+b_2-1}\sg_{b_1+b_2}\sg_{b_1+b_2+1} \sg_{b_1+b_2+2}} = 
1324$ ensures that $\sg_{b_1+2}$ is also less than the least element 
of $B_3$ and since $\sg_{b_1} < \sg_{b_1+2}$, we must have 
$\sg_{b_1} =x+2$ and $\sg_{b_1+2} = x+3$.  If we then remove 
the first $x+2$ cells which contain the numbers $1,\ldots, x+2$, 
then we must be left with a fixed point which has two bricks 
on $n-x-2$ cells. Then by our analysis of case 2, there 
are $n-x-2-3$ possibilities for $B_3$ so that we have 
a total of $\sum_{x=1}^{n-6} n-x-5 = \binom{n-5}{2}$ possibilities.

It follows that in subcase d.1 where $\tau = 1324$, 
we have $2(n-5) + \binom{n-5}{2}$ possibilities 
if $n \geq 7$ and no possibilities if $n < 7$. \\
\ \\
{\bf Subcase d.2.} $p \geq 5$. \\

In this case, we must have a $1324\ldots p$-match among the elements 
of bricks $B_1$ and $B_2$ which can only happen if 
$\red{\sg_{b_1-1}\sg_{b_1}\sg_{b_1+1}\dots \sg_{b_1+p-2}} = 
1324\ldots p$ and $b_2 \geq p-2$. 
Similarly we must have have $1324\ldots p$-match among the elements 
of bricks $B_2$ and $B_3$ which can only happen if 
$\red{\sg_{b_1+b_2-1}\sg_{b_1+b_2}\sg_{b_1+b_2+1}\dots \sg_{b_1+b_2+p-2}} = 
1324\ldots p$ which is a 
condition that does not involve $\sg_{b_1+1}$. Nevertheless, 
these conditions still force that $\sg_{b_1}, \sg_{b_1+1} <\sg_{b_1+b_2+1}$ so 
that if $\sg_{b_1-1} =x$, then $x$ is less than the least elements 
in bricks $B_2$ and $B_3$ so that $1, \ldots, x-1$ must be in 
brick $B_1$ and $\sg_{b_1+1} = x+1$.  We also have 
that $\sg_{b_1} < \sg_{b_1+2}$ and that $\sg_{b_2+2}$ must be less 
than the least element in brick $B_3$ which is $\sg_{b_1+b_2+1}$. 
It follows that it must be the case that 
$\sg_{b_1} =x+2$ and $\sg_{b_1+2} = x+3$.  If we then remove 
the first $x+p-3$ cells which contain the numbers $1,\ldots, x+p-3$, 
then we will be left with a fixed point which has two bricks 
on $n-x-p+3$ cells. Then by our analysis of case 2, there 
are $n-x-p+3-p+1$ possibilities for $B_3$ so that we have 
a total of $\sum_{x=1}^{n-2(p-2)-1} n-x-2p+4 = \binom{n-2p+4}{2}$ 
possibilities if $n \geq 2p-2$ and no possibilities if 
$n< 2p-2$.\\
\ \\

It follows that 
\begin{equation}
U_{1324,n}(y)|_{y^3} = \begin{cases} 0 & \mbox{if} \ n =0,1,2 \\
\binom{n-1}{2} & \mbox{if} \ n = 3,4 \ \mbox{and} \\
\binom{n-1}{2} + 2\binom{n-3}{2} + \binom{n-5}{2} + 2(n-5) = 2(n-3)^2 & \mbox{if} \ n \geq 5.
\end{cases}
\end{equation}
Similarly, if $p \geq 5$, then  
\begin{equation}
U_{1324\ldots p,n}(y)|_{y^3} = \begin{cases} 0 & \mbox{if} \ n =0,1,2 \\
\binom{n-1}{2} & \mbox{if} \ n = 3 \leq n \leq p \\
\binom{n-1}{2} + 2\binom{n-p+1}{2} & \mbox{if} \ n = p+1 \leq \ n \leq 2p-3
 \ \mbox{and} \\
\binom{n-1}{2} + 2\binom{n-p+1}{2} + \binom{n-2p+4}{2} & \mbox{if} \ n \geq 2p-2.
\end{cases}
\end{equation}

Note that $2(n-3)^2 = 2n(n-1) -10n +18$ so that 
\begin{eqnarray*}
U_{1324}(t,y)|_{y^3} &=&  \frac{t^3}{3!} + 3\frac{t^4}{4!} + 
\sum_{n \geq 5} 2(n-3)^2 \frac{t^n}{n!} \\
&=& \frac{t^3}{3!} + 3\frac{t^4}{4!} + 
 \sum_{n \geq 5} 2n(n-1)\frac{t^n}{n!} -\sum_{n \geq 5} 10n\frac{t^n}{n!} + 
+ \sum_{n \geq 5} 18\frac{t^n}{n!}  \\
&=& \frac{t^3}{3!} + 3\frac{t^4}{4!} + 
2t^2(e^t -\sum_{n=0}^2 \frac{t^n}{n!}) +
-10t(e^t -\sum_{n=0}^3 \frac{t^n}{n!})+18(e^t -\sum_{n=0}^4 \frac{t^n}{n!})\\
&=& (2t^2 -10t+18)e^t -18 -8t -t^2 +\frac{t^3}{3!} +\frac{t^4}{4!}.
\end{eqnarray*}

For $p \geq 5$, note that 
$$\binom{n-1}{2} + 2\binom{n-p+1}{2} + \binom{n-2p+4}{2} =
2n(n-1)+ (5-4p)n +3p^2-8p+7$$
so that 
\begin{multline*}
\sum_{n \geq 2p-2} \left(\binom{n-1}{2} + 2\binom{n-p+1}{2} + \binom{n-2p+4}{2}\right) 
\frac{t^n}{n!} = \\
\sum_{n \geq 2p-2} (2n(n-1)+ (5-4p)n +3p^2-8p+7) 
\frac{t^n}{n!} = \\
2t^2(e^t - \sum_{n=0}^{2p-5}\frac{t^n}{n!}) +  (5-4p)t(e^t - \sum_{n=0}^{2p-4}\frac{t^n}{n!}) + 
(3p^2-8p+7)(e^t - \sum_{n=0}^{2p-3}\frac{t^n}{n!}).
\end{multline*}
It follows that for $p \geq 5$, 
\begin{eqnarray*}
U_{1324\ldots p}(t,y)|_{y^3} &=& (2t^2+(5-4p)t+3p^2-8p+7)e^t + \\
&&\sum_{n=3}^p \binom{n-1}{2} \frac{t^n}{n!} + 
\sum_{n=p+1}^{2p-3} (\binom{n-1}{2} + 2\binom{n-p+1}{2})\frac{t^n}{n!} - \\
&&\left(2t^2  \sum_{n=0}^{2p-5}\frac{t^n}{n!} + (5-4p)t  \sum_{n=0}^{2p-4}\frac{t^n}{n!} +(3p^2-8p+7)  \sum_{n=0}^{2p-3}\frac{t^n}{n!}\right).
\end{eqnarray*}
One can then use Mathematica to show that 
\begin{equation}
U_{1324\ldots p}(t,y)|_{y^3} = (2t^2-(5-4p)t+3p^2-8p+7)e^t + \sum_{n=0}^{2p-3} 
f(n,p)
\end{equation}
where 
\begin{equation}\label{fnpdef}
f(n,p)= \begin{cases} -3p^2+8p-7  \ \mbox{if} \ n=0, \\
 -3p^2+12p-12 \ \mbox{if} \ n=1, \\
 -3p^2+16p-21 \ \mbox{if} \ n=2, \\
-\frac{3n^2}{2}+(4p-\frac{9}{2}) -3p^2 +8p-6  \ \mbox{if} \ 3 \leq n \leq p, 
\ \mbox{and} \\
-\frac{n^2}{2}+(2p-\frac{7}{2}) -3p^2 +7p-6 \  
\mbox{if} \ p+1 \leq n \leq 2p-3.
\end{cases}
\end{equation}

Now for any $\tau$, we can 
write 
\begin{eqnarray} 
NM_{\tau}(t,1,y) &=& \frac{1}{U_\tau(t,y)} = 
\frac{1}{1 -(A_\tau(t)y -B_\tau(t)y^2 + C_\tau(t)y^3 +O(y^4))}  \nonumber \\
&=& 1+\sum_{n \geq 1} (A_\tau(t)y -B_\tau(t)y^2 + C_\tau(t)y^3 +O(y^4))^n.
\end{eqnarray}
It then follows that 
\begin{eqnarray*}
NM_{\tau}(t,1,y)|_y &=& A_\tau(t), \\
NM_{\tau}(t,1,y)|_{y^2} &=& (A_\tau(t))^2 - B_\tau(t), \ \mbox{and} \\
NM_{\tau}(t,1,y)|_{y^3} &=& (A_\tau(t))^2 - 2A_\tau(t)B_\tau(t) + C_\tau(t).
\end{eqnarray*}

We have shown that 
\begin{eqnarray*}
A_{1324}(t) &=& e^t-1, \\
B_{1324}(t) &=& (2t-4)e^t+4+2t+\frac{t^2}{2}, \ \mbox{and}\\
C_{1324}(t) &=& (2t^2 -10t+18)e^t -18 -8t -t^2 +\frac{t^3}{3!} +\frac{t^4}{4!}.
\end{eqnarray*}
One can then use Mathematica to compute that 
\begin{equation}
D^{(1)}_4(t) = NM_{1324}(t,1,y)|_{y^2} = e^{2t}-(2t-2)e^t-3 -2t +\frac{t^3}{6}.
\end{equation}
It follows that for $n \geq 4$, 
\begin{equation}\label{1324-1d}
d^{(1)}_{n,4}  = 2^n-2n+2.
\end{equation}
This is easy to explain directly. That is, if $\sg \in NM_{1324,n}$ 
and has one descent, then $\sg$ has either one or two left-right-minima. 
Thus for $p \geq 4$, 
\begin{eqnarray*}
d^{(1)}_{n,4} &=& NM_{1324,n}(x,y)|_{xy^2} + NM_{1324,n}(x,y)|_{x^2y^2} \\
&=& 2^{n-1}-2n +4 -1 + 2^{n-1}-1\\
&=& 2^n-2n+2.
\end{eqnarray*}

One can also use Mathematica or Maple to compute that 
\begin{equation}
D^{(2)}_4(t) = NM_{1324}(t,1,y)|_{y^3} = e^{3t} +(5-4t)e^{2t}+(t^2-10t+5)e^t-11 -4t +\frac{t^3}{6} +\frac{t^4}{24}.
\end{equation}
It then follows that for $n \geq 5$, 
\begin{equation}\label{1324-2d}
d^{(2)}_{n,4} = 
3^n +(5-2n)2^n +n^2-11n+5.
\end{equation}
We do not know of a simple direct proof of this result.

We have shown that for $p \geq 5$, 
\begin{eqnarray*}
A_{1324\ldots p}(t) &=& e^t-1, \\
B_{1324\ldots p}(t)  &=& (2t-p)e^t+p+\sum_{n=1}^{p-2}(p-n-1) \frac{t^n}{n!}, \ \mbox{and}\\
C_{1324\ldots p}(t) &=& (2t^2 +(5-4p)t+3p^2-8p+7)e^t +\sum_{n=0}^{2p-3} f(n,p) \frac{t^n}{n!}.
\end{eqnarray*}
One can then use Mathematica to compute that 
\begin{equation}
D^{(1)}_p(t) = NM_{1324}(t,1,y)|_{y^2} = e^{2t}-(2t-p+2)e^t - \sum_{n=0}^{p-2} (p-1-n)\frac{t^n}{n!}.
\end{equation}
It follows that for $n \geq p-1$, 
\begin{equation}\label{1324p-1d}
d^{(1)}_{n,p} = 2^n-2n-p+2.
\end{equation}
One can easily modify the direct argument that we used to prove 
$d^{(1)}_{n,4} = 2^n-2n+2$ for $n \geq 4$ to give a direct 
proof of this result.

One can also use Mathematica to compute that 
\begin{multline}
D^{(2)}_p(t) = NM_{1324\ldots p}(t,1,y)|_{y^3} = 
e^{3t} +(2p-3-4t)e^{2t}+\\
\left(3p^2-12p+10 +(13-6p)t + (5-p)t^2 - \sum_{n=3}^{p-2} (2(p-n-1)\frac{t^n}{n!} \right)e^t+\\
2p -1+ \sum_{n=3}^{p-2} (2(p-n-1)\frac{t^n}{n!} + \sum_{n=0}^{2p-3} f(n,p) \frac{t^n}{n!}
\end{multline}
where $f(n,p)$ is defined as in (\ref{fnpdef}). 
For example, one can compute that 
\begin{eqnarray*}
D^{(2)}_5(t) &=&  e^{3t} +(7-4t)e^{2t} + (25 -17t -\frac{2t^3}{3!})
e^t -33 -21t -6t^2 -t^3 -\frac{3t^4}{4!}-\frac{t^5}{5!}, \\
D^{(2)}_6(t) &=&  e^{3t} +(9-4t)e^{2t} + 
(46 -23t -t^2 - \frac{4t^3}{3!} -\frac{2t^4}{4!})e^t \\
&&-56 -40t -\frac{27t^2}{2} -\frac{17 t^3}{3!} -\frac{3t^4}{4!}-
\frac{6t^5}{5!} -\frac{3t^6}{6!} -\frac{t^7}{7!}, \ \mbox{and}\\
D^{(2)}_7(t) &=&  e^{3t} +(11-4t)e^{2t} + (73 -29t -2t^2 - t^3 - 
\frac{4t^4}{4!} -\frac{2t^5}{5!})e^t \\
&&-85 -65t -\frac{48t^2}{2}  -\frac{34t^3}{3!} -\frac{23t^4}{4!}-
\frac{15t^5}{5!} -\frac{10t^6}{6!} -\frac{6t^7}{7!} -\frac{3t^8}{8!} 
-\frac{t^9}{9!}.
\end{eqnarray*}
It then follows that for $n \geq 2p-2$, 
\begin{multline}\label{1324p-2d}
d^{(2)}_{n,p} = 3^n +(2p-3-2n)2^n +\\
3p^2 -12p +10 +(13-6p)n + (5-p)n(n-1) - 
\sum_{k=3}^{p-2} \frac{2(p-k-1)}{k!}n(n-1) \cdots (n-k+1).
\end{multline}
For example, for $n \geq 8$, 
$$d^{(2)}_{n,5} 
= 3^n+(7-2n)2^n + 25 -\frac{53n}{3} +n^2 - \frac{n^3}{3}.$$
For $n \geq 10$, 
$$d^{(2)}_{n,6} 
= 3^n+(9-2n)2^n + 46 -\frac{136n}{6} +\frac{n^2}{12} - \frac{n^3}{6} -\frac{n^4}{12}.$$
For $n \geq 12$, 
$$d^{(2)}_{n,7} 
= 3^n+(11-2n)2^n + 73 -\frac{142n}{5} +\frac{7n^3}{12} - \frac{n^5}{60}.$$

\section{Conclusions}

In this paper, we showed that if $\tau$ is a permutation 
which starts with 1, then 
$NM_{\tau}(t,x,y)$ is always of the form 
$\left( \frac{1}{U_\tau(t,y)}\right)^x$ where 
$U_\tau(y) = 1 + \sum_{n \geq 1} U_{\tau,n}(y) \frac{t^n}{n!}$. 
In the special case where $\tau$ has one descent, we 
showed how to use the homomorphism method to give a  
simple combinatorial description of $U_{\tau,n}(y)$ for 
any $n \geq 1$. We then used 
this combinatorial description to show that 
the $U_{\tau,n}(y)$s satisfy simple recursions in the 
special case where $\tau = 1324 \ldots p$ and $p \geq 4$.  

The methods introduced in this paper can be use to prove 
several other similar results for other collections 
of patterns that start with 1 and have one descent.
For example, suppose that $\tau=1p2\dots(p-1)$ where 
$p \geq 4$. Then we can prove that 
$$NM_\tau(t,x,y)=\left(\frac1{U_\tau(t,y)}\right)^x \text{ where }U_\tau(t,y)=1+\sum_{n\geq1}U_{\tau,n}(y)\frac{t^n}{n!},$$
$U_{\tau,1}(y)=-y$, and for $n\geq 2$,
$$U_{\tau,n}(y)=(1-y)U_{\tau,n-1}(y)+\sum_{k=1}^{\lfloor\frac{n-2}{p-2}\rfloor}(-y)^{k}{n-k(p-3)-2 \choose k-1}U_{\tau,n-(k(p-2)+1)}(y).$$
If $p \geq 5$ and $\tau = 134 \ldots (p-1)2p$, then we can 
prove that 
$$NM_\tau(t,x,y)=\left(\frac1{U_\tau(t,y)}\right)^x \text{ where }U_\tau(t,y)=1+\sum_{n\geq1}U_{\tau,n}(y)\frac{t^n}{n!},$$
$U_{\tau,1}(y)=-y$, and for $n\geq 2$,
$$U_{\tau,n}(y)=(1-y)U_{\tau,n-1}(y)+\sum_{k=1}^{\lfloor \frac{n-2}{p-2} \rfloor}(-y)^k  
\frac{1}{(p-3)k+1}\binom{k(p-2)}{k} U_{\tau,n-k(p-2)-1}.$$
These results will appear in subsequent papers.

\end{document}